\documentclass[12pt]
{amsart}
\usepackage{color}

\usepackage{verbatim}
\usepackage{amsmath}
\usepackage{amsthm}
\usepackage{amssymb}
\usepackage{enumerate}

\newif\ifdeveloping

\let\QED\qed
\newcommand{\prlabel}[1]{\renewcommand{\qed}{\QED${}_{\ref{#1}}$}}
\newcommand{\prtxtlabel}[1]{\renewcommand{\qed}{\QED${}_{{\mbox{\tiny #1}}}$}}
\newcommand{\prnolabel}{\prtxtlabel{}}

\newtheorem{theorem}{Theorem}[section]
\newtheorem{mtheorem}[theorem]{Main theorem}
\newtheorem{proposition}[theorem]{Proposition}
\newtheorem{lemma}[theorem]{Lemma}
\newtheorem{corollary}[theorem]{Corollary}
\newtheorem{fact}[theorem]{Fact}

\newtheorem{claim}{Claim}[theorem]
\newtheorem{problem}[theorem]{Problem}
\newtheorem{case}{Case}

\newtheorem{obs}[theorem]{Observation}
\newtheorem{qtheorem}{Theorem}

\theoremstyle{definition}

\newtheorem{definition}[theorem]{Definition}
\newtheorem{mdefinition}[theorem]{Main Definition}

\theoremstyle{remark}
\newtheorem{remark}{Remark}  
{}
\newcommand{\eqs}{=^*}
\newcommand{\subss}{\subs^*}

\ifdeveloping
\usepackage[notref,notcite]{showkeys}
\fi

\newcommand{\prtime}{{\count0=\time\divide\count0 by 60
\count1=-\count0\multiply\count1 by 60
\advance\count1 by \time
\the\count0:\the\count1}
}

\def\myheads#1;#2;{
\pagestyle{myheadings}
\markboth{{\sc\hfill #1\hfill\protect\makebox[0cm][r]{\rm\today; \prtime}}}
{{\sc\protect\makebox[0cm][l]{\rm\today;\ \prtime}\hfill #2\hfill}}
\thispagestyle{myheadings}
}

\newcommand{\acal}{{\mathcal A}}

\newcommand{\ccal}{{\mathcal C}}
\newcommand{\dcal}{{\mathcal D}}
\newcommand{\ecal}{{\mathcal E}}
\newcommand{\fcal}{{\mathcal F}}
\newcommand{\ical}{{\mathcal I}}
\newcommand{\ncal}{{\mathcal N}}
\newcommand{\pcal}{{\mathcal P}}
\newcommand{\tcal}{{\mathcal T}}
\newcommand{\ucal}{{\mathcal U}}
\newcommand{\vcal}{{\mathcal V}}

\newcommand{\setm}{\setminus}
\newcommand{\empt}{\emptyset}
\newcommand{\subs}{\subset}

\newcommand{\oo}{{{\omega}_1}}
\newcommand{\dom}{\operatorname{dom}}
\newcommand{\nden}{\operatorname{\ncal}}

\def\<{\left\langle}
\def\>{\right\rangle}
\def\cf{\operatorname{cf}}
\def\br#1;#2;{\bigl[ {#1} \bigr]^ {#2} }
\newcommand{\fn}{\operatorname{Fn}}

\newcommand{\newcases}{\setcounter{case}{0}}

\newcommand{\restr}
{\mathop{\hspace{0.01ex}|\hspace*{-0.02ex}{\grave{}}\hspace{0.4ex}}}
\theoremstyle{plain}

\def\separ{separating}
\def\finc{{\mathbb{F}\mathbb{I}\makebox[-1pt]{}\mathbb{N}}}
\def\celh{\operatorname{\hat c}}
\def\cel{\operatorname{c}}
\newcommand{\conc}{^\frown}

\newcommand{\ores}{open hereditarily irresolvable}
\newcommand{\hres}{hereditarily resolvable}
\newcommand{\res}{resolvable}
\newcommand{\hires}{hereditarily irresolvable}

\newcommand{\eres}{extraresolvable}

\newcommand{\seres}{strongly extraresolvable}
\newcommand{\seresk}[1]{strongly $#1$-extraresolvable}

\newcommand{\hyres}{hyperresolvable}
\newcommand{\whyres}{fragmented}
\newcommand{\nwd}{\operatorname{nwd}}
\newcommand{\nwds}{\operatorname{\ncal}}

\newcommand{\appr}{limit}
\newcommand{\bbb}{\mathbb B}
\newcommand{\cbb}{\mathbb C}
\newcommand{\dbb}{\mathbb{D}}
\newcommand{\ebb}{\mathbb{E}}
\newcommand{\fbb}{\mathbb{F}}
\newcommand{\dis}{\operatorname{D}}
\newcommand{\piweight}{\operatorname{\pi w}}

\newcommand{\mos}[2]{\operatorname{\frak {M
}}(#1,#2)}

\newcommand{\mosaicu}[2]{$(#1,#2)$-mosaic}
\newcommand{\mosaic}[1]{$#1$-mosaic}
\newcommand{\pieu}[2]{$(#1,#2)$-piece}
\newcommand{\pie}[1]{$#1$-piece}
\newcommand{\pis}[2]{\pcal(#1,#2)}
\newcommand{\pisx}[1]{\pcal(#1)}
\newcommand{\irrac}{\mathbb P}

\newcommand{\hint}[1]{\operatorname{\ical}(#1)}

\begin{document}

\author[I. Juh\'asz]{Istv\'an Juh\'asz}
\address{Alfr{\'e}d R{\'e}nyi Institute of Mathematics}
\email{juhasz@renyi.hu}

\author[L. Soukup]{
Lajos Soukup }

\thanks{The preparation of this paper was supported by the
Hungarian National Foundation for Scientific Research grant no. 37758}
\thanks{
The second author was also partially
supported by Grant-in-Aid for JSPS Fellows No.\ 98259 of the Ministry of
Education, Science, Sports and Culture, Japan and 
by the Bolyai Fellowship}
\address{
Alfr{\'e}d R{\'e}nyi Institute of Mathematics }
\email{soukup@renyi.hu}

\author[Z. Szentmikl\'ossy]{
Zolt\'an Szentmikl\'ossy}
\email{zoli@renyi.hu}

\title[$\dcal$-forced spaces]{$\dcal$-forced spaces: a new approach to resolvability}

\begin{abstract}
We introduce a ZFC method that enables us to build
spaces (in fact special dense subspaces of certain Cantor cubes) 
in which we have "full control" over all dense
subsets. 

Using this method we are able to construct, in ZFC, 
for each uncountable regular cardinal $\lambda$ 
a $0$-dimensional $T_2$, hence Tychonov, space which is $\mu$-resolvable
for all $\mu < \lambda$ but not
$\lambda$-resolvable. This yields the final (negative) solution of a celebrated 
problem of Ceder and Pearson raised in 1967: Are $\omega$-resolvable
spaces maximally resolvable?
This method enables us to solve several other open problems concerning
resolvability as well.
\end{abstract}

\maketitle
\ifdeveloping
\myheads{Resolvable};{Resolvable};
\fi

\section{Introduction}

Resolvability questions about topological spaces were first
studied by E. Hewitt, \cite{He}, in 1943. Given a cardinal
${\kappa}>1$, a topological space $X=\<X,{\tau}_X\>$ is called
{\em ${\kappa}$-resolvable} iff it contains ${\kappa}$ disjoint
dense subsets. $X$ is {\em resolvable} iff it is 2-resolvable and
{\em irresolvable} otherwise.

If $X$ is ${\kappa}$-resolvable and $G\subs X$ is any non-empty
open set in $X$ then clearly ${\kappa}\le |G|$.
Hence if $X$ is ${\kappa}$-resolvable then we have
$\kappa \le \operatorname{\Delta}(X)$ where
$$\operatorname{\Delta}(X) =
\min\bigl\{|G|:G\in {\tau}_X\setm \{\empt\}\bigr\}.$$
This observation explains the
following  terminology of J.Ceder, \cite{Ce}:
a space  $X$ is called {\em maximally resolvable } iff it is
$\operatorname{\Delta}(X)$-resolvable.

Ceder and Pearson, in \cite{CP}, raised the question whether an
${\omega}$-resolvable space is necessarily  maximally resolvable?
El'kin,  \cite{El}, Malykhin, \cite{Ma2}, Eckertson, \cite{E}, and
Hu,  \cite{Hu}, gave several counterexamples, but either these
spaces
were not  even $T_2$ 
or their construction was not carried out in ZFC. Our theorems
\ref{tm:mplu} and \ref{tm:limit} give a large number of
  $0$-dimensional $T_2$ (and so  Tychonov) counterexamples in ZFC.
The question if this can be done has been asked much more recently
again in \cite{Co} and \cite{CoHu}.

Our results are obtained with the help of a new method that is
presented
in section \ref{sc:main}. Here we first introduce the new and
simple concept of {\em $\dcal$-forced spaces}.
Given a family $\dcal$ of dense subsets of the space $X$ we say
that   $X$ (or its topology) is {\em $\dcal$-forced} if any subset
of $X$ can only be dense in $X$ if
$\dcal$ forces this to happen. The exact formulation of this reads
as follows: If $S$ is dense in $X$ then $S$ includes a set of the
form
$$M=\cup\{V\cap D_V:V\in\vcal\}$$ where $\vcal$ is a maximal
disjoint collection of open sets in $X$ and $D_V \in \dcal$ for
all $V \in \vcal$. Such a set $M$, that is clearly dense in $X$,
is called a $\dcal$-mosaic. Then, in lemmas \ref{i:nwd},
\ref{i:irres}, \ref{i:space}, and 2.10, we establish the basic
properties of $\dcal$-forced spaces.

In the next section we prove our main result, theorem
\ref{tm:main},  that will allow us to construct $\dcal$-forced
subspaces of certain Cantor cubes with a wide range of
resolvability, resp. irresolvability properties. Thus, in sections
4 and 5, we shall be able to answer not only the problem of Ceder
and Pearson mentioned above but several other open problems as well, like
\cite[Question 4.4]{ASTTW}, \cite[Problem 8.6]{AC},
\cite[Questions 3.4, 3.6, 4.5]{E}, or a problem of Comfort and Hu
mentioned in \cite[Discussion 1.4]{CoHu}.

In the remaining part of this introduction we summarize our
further notation and terminology, most of it is standard.
 
A space $X$ is called {\em \ores\ (OHI) } iff every nonempty open
subspace of $X$ is irresolvable. It is well-known that every
irresolvable space has a non-empty open subspace that is OHI.
Clearly, $X$ is OHI iff every dense subset of $X$ contains a dense
open subset, i. e. if $S\subs X$ dense in $X$ implies that
$\operatorname{Int}(S)$ is dense, as well.

Next, a space $X$ is called {\em \hires (HI)} iff all subspaces of
$X$ are irresolvable. Since a space having an isolated point is
trivially irresolvable, any space is HI iff all its {\em crowded}
subspaces are irresolvable. (Following van Douwen, we call a space
crowded if it has no isolated points.) Having this in mind, if P
is any resolvability or irresolvability property of topological
spaces then the space $X$ is called {\em hereditarily P} iff all
{\em crowded} subspaces of $X$ have property P.

Following the terminology of  \cite{vD2}, a topological space $X$
is called {\em NODEC} if all nowhere dense subsets of $X$ are
closed, and hence closed discrete. All spaces obtained by our main
theorem \ref{tm:main} will be NODEC.

A space is called {\em submaximal } (see \cite{He}) iff all of its
dense subsets are open. The following observation is easy to prove
and will be used repeatedly later: a space is submaximal iff it is
both OHI and NODEC.

A set $D\subs X$ is said to be {\em ${\kappa}$-dense} in $X$ iff
$|D\cap U|\ge {\kappa}$ for each nonempty open set $U\subs X$.
Thus $D$ is dense iff it is 1-dense. Also, it is obvious that the
existence of a $\kappa$-dense set in $X$ implies
$\operatorname{\Delta}(X) \ge \kappa$.

We shall denote by $\nden(X)$ the family of all nowhere dense
subsets of a space $X$. Clearly, $\nden(X)$ is an ideal of subsets
of $X$ and the notation $=^*$ or $\subset^*$ will always be used
to denote equality, resp. inclusion modulo this ideal.

Following the notation introduced in \cite{CG2}, we shall write
$$\nwd(X)=\min\{|Y|:Y\in \pcal(X)\setm \nden(X)\} = non-(\nden(X)),$$ i. e.
$\nwd(X)$ is the minimum cardinality of a
 somewhere dense subset of $X$.

Malychin was the first to suggest studying families of dense sets
of a space $X$ that are  {\em almost disjoint} with respect to the
ideal $\nden(X)$ rather than disjoint, see \cite{Ma3}. He calls a
space $X$ {\em\eres} if there are $\operatorname{\Delta}(X)^+$
many dense sets in $X$ such that any two of them have nowhere
dense intersection. Here we generalize this concept by defining a
space $X$ to be $\kappa$-{\em extraresolvable} if there are
$\kappa$ many dense sets in $X$ such that any two of them have
nowhere dense intersection. (Perhaps $\kappa$-almost resolvable
would be a better name for this.) Note that, although
$\kappa$-extraresolvability of $X$ is mainly of interest if
$\kappa > \operatorname{\Delta}(X)$, it does make sense for
$\kappa \le \operatorname{\Delta}(X)$ as well.  Clearly,
$\kappa$-resolvable implies $\kappa$-extraresolvable, moreover the
converse holds if $\kappa = \omega$, however we could not decide
if these two concepts coincide if $$\omega < \kappa \le
\operatorname{\Delta}(X).$$ In particular, we would like to know
the answer to the following question.
\begin{problem}
Let $X$ be an extraresolvable ($T_2$, $T_3,$ or Tychonov) space with
$\operatorname{\Delta}(X) \ge \omega_1$. Is $X$ then $\omega_1$-
resolvable?
\end{problem}
Note that a counterexample to Problem 1.1 is also a counterexample
to the Ceder-Pearson problem.

Finally we  mention a variation of extraresolvability. 
The space $X$ is called {\em \seresk{\kappa}}
iff there are $\kappa$ many dense subsets
$\{D_{\alpha}:{\alpha}<{\kappa}\}$ of $X$ such that
$|D_{\alpha}\cap D_{\beta}|<\nwd(X)$  whenever
$\{{\alpha},{\beta}\}\in \br {\kappa};2;$. We say that $X$ is
{\em\seres} iff it is \seresk{\operatorname{\Delta}(X)^+}.
Clearly, \seresk{(\kappa)} implies $(\kappa)$- extraresolvable.

\section{$\dcal$-forced spaces}\label{sc:main}

\begin{definition}
Let $\dcal$ be a family of dense subsets of a space $X$.  A subset
$M\subs X$ is called a {\em \mosaicu {\dcal}X} iff
there is a maximal disjoint family
$\vcal$ of open subsets of $X$ and for each $V\in\vcal$
there is $D_V\in \dcal$
such that
\begin{displaymath}
M=\cup\{V\cap D_V:V\in\vcal\}.
\end{displaymath}
A set $M$ of the above form with $\vcal$ disjoint, but not
necessarily maximal disjoint, is called a {\em partial} {\mosaicu {\dcal}X}.

A set $P$ of the form $P = D\cap U$, where $D\in \dcal$ and $U$ is
a nonempty open subset of $X$, is called a {\em \pieu{\dcal}X}.
So, naturally, any (partial) { \mosaicu {\dcal}X} is composed of
$(\dcal,X)$-pieces. Let
\begin{displaymath}
\mos {\dcal}{X}=\{M:\text{ $M$ is a \mosaicu {\dcal}X}\}
\end{displaymath}
and
\begin{displaymath}
\pis {\dcal}X=\{P:\text{$P$ is a \pieu{\dcal}X}\}.
\end{displaymath}
\end{definition}

When the space $X$ is clear from the context
we will omit it from the notation: we will write  {\em \mosaic{\dcal}}
 instead of  \mosaicu{\dcal}X, and  {\em \pie{\dcal}} instead of
\pieu{\dcal}X, etc. The following statement is now obvious.

\begin{fact}
Every \mosaicu {\dcal}X is dense in $X$ and every \pieu{\dcal}X is
somewhere dense in $X$.
\end{fact}

Thus we arrive at the following very simple but, as it turns out,
very useful concept.

\begin{mdefinition}
Let $\dcal$ be a family of dense subsets of a topological space
$X$. We say  that the space $X$ (or its topology)
 is  {\em ${\dcal}$-forced} iff
every dense subset $S$ of $X$ includes a $\dcal$-mosaic $M$, i. e.
there is $M \in \mos {\dcal}{X}$ with $M \subset S$ .
\end{mdefinition}

It is easy to check that one can give the following alternative
characterization of being ${\dcal}$-forced.

\begin{fact}\label{f:pd}
The space $X$ is $\dcal$-forced iff every somewhere dense subset
of $X$ includes a \pieu{\dcal}X.
\end{fact}

Since $X$ is always dense in $X$, the simplest choice for $\dcal$
is $\{X\}$.

\begin{fact}
 A subset $P\subs X$ is an \pie{\{X\}}
 iff it is non-empty  open; $M$ is an \mosaic{\{X\}} iff it is dense open in $X$.
Consequently, $X$ is $\{X\}$-forced iff it is OHI.
\end{fact}

Let us now consider a few further, somewhat less obvious,
properties of $\dcal$-forced spaces. The first result yields a
useful characterization of nowhere dense subsets in such spaces.
Note that a subset $Y$ of any space $X$ is nowhere dense iff $S
\backslash Y$ is dense in $X$ for all dense subsets $S$ of $X$.
Not surprisingly, in a $\dcal$-forced space it suffices to check
this for members of $\dcal$.
\begin{lemma}\label{i:nwd}
Assume that $X$ is $\dcal$-forced. Then
  \begin{displaymath}
  \nwds(X)=\{Y\subs X: D\setm Y
      \text{ is dense in $X$ for each $D\in \dcal$}\}.
  \end{displaymath}
\end{lemma}

\begin{proof}
Assume that $Y\notin \nwds(X)$, i.e. $Y$ is somewhere dense. Then,
by fact \ref{f:pd}, $Y$ contains some \pie{\dcal}  $U\cap D$,
where $D\in\dcal$ and $U$ is a nonempty open subset of $X$. Then
$(D\setm Y)\cap U=\empt$, i.e. $D\setm Y$ is not dense. This
proves that the right-hand side of the equality includes the left
one. The converse inclusion is obvious.
\end{proof}

The following result will be used to produce irreducible (even
OHI) spaces. Of course, the superscript * in its formulation
designates equality and inclusion modulo the ideal $\nwds(X)$ of
nowhere dense sets.

\begin{lemma}\label{i:irres}
Let $X$ be $\dcal$-forced and $S\subs X$ be dense such that
  \begin{equation}\tag{$\dag$}
 \text{
for each $D\in \dcal$ we have $S\cap D\eqs \empt$  or $S\subss D$.
}
  \end{equation}
 Then $S$, as a subspace of $X$, is OHI.
\end{lemma}

\begin{proof}
Let $T\subs S$ be dense in $S$, then $T$ is also dense in $X$,
hence it must contain a \mosaic{\dcal}, say $M=\bigcup\{V\cap
D_V:V\in \vcal\}$. But then we have $S\subs^* D_V$ for each
$V\in\vcal$ by $(\dag)$. Consequently, $$T\cap V\subs S\cap
V\subss V\cap D_V\subs T\cap V$$ and so
 $T\cap V=^*S\cap V$ holds for all $V\in \vcal$. This clearly implies that $T =^*
 S$. In other words, we have shown that every dense subset $T$ of $S$
has nowhere dense complement in $S$,
i. e. the subspace $S$ of $X$ is OHI.
\end{proof}

The following lemma will enable us to conclude that certain
$\dcal$-forced spaces are not $\kappa$-(extra)resolvable for
appropriate cardinals $\kappa$.

\begin{lemma}\label{i:space}
Assume that $X$ is a topological space and $\dcal$ is a family of
dense subsets of $X$. Assume, moreover, that ${\mu}\ge \celh(X)$
(i.e. $X$ does not contain ${\mu}$ many pairwise disjoint open
subsets) and
\begin{multline}
\tag{$*$}
\text{for each $\ecal\in \br \dcal;{\mu};$ there is
$\fcal\in \br \ecal;\celh(X);$
such that}\\
\text{ $D_0\cap D_1$  is dense in $X$ whenever $\{D_0,D_1\}\in \br
\fcal;2;$}.
\end{multline}
Then for any family of $\dcal$-pieces $\{P_i:i<{\mu}\}\subs
\pisx{\dcal}$ there is $\{i,j\}\in \br {\mu};2;$ such that
$P_i\cap P_j$ is somewhere dense in $X$.

In particular, if  $X$ is
 $\dcal$-forced and $|\dcal|^+\ge \celh(X)$ then $X$ is not
$|\dcal|^+$-extraresolvable (hence not $|\dcal|^+$-resolvable,
either).
\end{lemma}

\begin{proof}
  Assume that $P_i=U_i\cap D_i$, where $D_i\in\dcal$ and $U_i$ is a nonempty
open  subset of $X$ for all $i \in \mu$. By $(*)$ there is $I\in
\br \mu;\celh(X);$ such that $D_i\cap D_j$ is dense for each
$\{i,j\}\in \br I;2;$. By the definition of $\celh(X)$, there is
$\{i,j\}\in \br I;2;$ such that $U= U_i\cap U_j$ is non-empty. But
then $U\cap D_i\cap D_j\subs P_i\cap P_j$, hence $P_i\cap P_j$ is
dense in the nonempty open set $U$.

The last statement now follows because $\dcal$ trivially satisfies
condition $(*)$ with $\mu = |\dcal|^+$ and, as $X$ is
$\dcal$-forced, every dense subset of $X$ includes a $\dcal$-piece
(even a $\dcal$-mosaic).
\end{proof}

The following fact is obvious.

\begin{fact}\label{i:res}
Let $\dcal$ be a family of dense sets in $X$ and
$$M=\bigcup\{V\cap D_V:V\in \vcal\}$$ be a partial $\dcal$-mosaic.
If all the dense sets $D_V$ are ${\mu}$-(extra)resolvable for $V
\in \vcal $ then so is $M$.
\end{fact}
We finish this section with a result that, together with fact
\ref{i:res}, will be used to establish hereditary
(extra)resolvability properties of several examples constructed
later.
\begin{lemma}
Let $X$ be a $\dcal$-forced space in which every crowded subspace
is somewhere dense. (This holds e. g. if $X$ is NODEC.) Then for
every crowded $S \subset X$ there is a partial $\dcal$-mosaic $M
\subset S$ that is dense in $S$. So if, in addition, all $D \in
\dcal$ are ${\mu}$-resolvable (resp. ${\mu}$-extraresolvable) then
$X$ is hereditarily ${\mu}$-resolvable (resp.
${\mu}$-extraresolvable).
\end{lemma}

\begin{proof}
Let $\vcal$ be a maximal disjoint family of open sets $V$ such
that there is an element $D_V \in \dcal$ with $V \cap D_V \subset
S$ and consider the partial $\dcal$-mosaic
$$M=\bigcup\{V\cap D_V:V\in \vcal\}.$$ Then $M \subset S$ is dense in $S$, since
otherwise, in view of the maximality of $\vcal$, the set $S
\setminus \overline{M} \neq \emptyset$ would be crowded and could
not include any $\dcal$-piece. The last sentence now immediately
follows using fact \ref{i:res}.
\end{proof}

\section{The Main Theorem}

We have introduced the concept of $\dcal$-forced spaces but one
question that immediately will be raised is if there are any
beyond the obvious choice of $\dcal = \{X\}$? The aim of this
section is to prove theorem \ref{tm:main} that provides us with a
large supply of such spaces. All these spaces will be dense
subspaces of Cantor cubes, i. e. powers of the discrete two-point
space $D(2)$. As is well-known, there is a natural one-to-one
correspondence between dense subspaces of size $\kappa$ of the
Cantor cube $D(2)^\lambda$ and independent families of
2-partitions of $\kappa$ indexed by $\lambda$. (A partition of a
set $S$ is called a {\em ${\mu}$-partition} if it partitions $S$
into ${\mu}$ many pieces.) For technical reasons, we shall produce
our spaces by using partitions rather than Cantor cubes.

We start with fixing some notation and terminology.

Let $\vec {\lambda}=\<{\lambda}_{\zeta}:{\zeta}<{\mu}\>$ be a
sequence of cardinals. We set
\begin{multline}\notag
\finc(\vec {\lambda})=\{{\varepsilon}:\text{${\varepsilon}$ is a
finite function with }
\dom {\varepsilon}\in \br {\mu};<{\omega};\text{ and }\\
{\varepsilon}({\zeta})\in {\lambda}_{\zeta} \text{ for all
${\zeta}\in\dom {\varepsilon}$}\}.
\end{multline}
Note that if $\lambda_\zeta = \lambda$ for all $\zeta < \mu$ then
$$\finc(\vec {\lambda}) = Fn(\mu,\lambda).$$

Let $S$ be  a set, $\vec
{\lambda}=\<{\lambda}_{\zeta}:{\zeta}<{\mu}\>$ be a  sequence of
cardinals, and
$\bbb=\left\{\<B_{\zeta}^i:i<{\lambda}_{\zeta}\>:{\zeta}<{\mu}
\right\}$ be a family of partitions of $S$.
Given a cardinal ${\kappa}$
we say that
$\bbb$ is {\em ${\kappa}$-independent} iff
\begin{displaymath}
\bbb[\varepsilon]\stackrel{def}=
\bigcap\{B^{{\varepsilon}({\zeta})}_{\zeta}:{\zeta}\in\dom
{\varepsilon} \}
\end{displaymath}
has cardinality at least ${\kappa}$
for each ${\varepsilon}\in \finc(\vec {\lambda})$. 
$\bbb$ is {\em independent}
iff it is 1-independent, i.e. $\bbb[{\varepsilon}]\ne \empt$
for each ${\varepsilon}\in \finc(\vec {\lambda})$. 
$\bbb$ is {\em\separ} iff for each $\{{\alpha},{\beta}\}\in \br S
;2;$ there are ${\zeta}<{\mu}$ and $\{{\rho},{\nu}\}\in
\br{\lambda}_{\zeta};2;$ such that ${\alpha}\in B^{\rho}_{\zeta}$
and ${\beta}\in B^{\nu}_{\zeta}$.

We shall denote by $\tau_\bbb$ the (obviously zero-dimensional)
topology on $S$ generated by the subbase $\{B_{\zeta}^i: \zeta <
\mu,\ i < {\lambda}_{\zeta}\}$, moreover we set
$X_{\bbb}=\<{S},{\tau}_{\bbb}\>$. Clearly, the family
$\{\bbb[{\varepsilon}]:{\varepsilon}\in \finc(\vec {\lambda})\}$
is a base for the space $X_{\bbb}$. Note that $X_{\bbb}$ is
Hausdorff iff $\bbb$ is separating.

The following statement is very easy to prove and is well-known.
It can certainly be viewed as part of the folklore.

\begin{obs}\label{lm:embed}
Let ${\kappa}$ and ${\lambda}$ be infinite cardinals. Then, up to
homeomorphisms, there is a natural one-to-one correspondence
between dense subspaces $X$ of ${D(2)}^{\lambda}$ of size
${\kappa}$ and  spaces of the form $X_{\bbb} =
\<{\kappa},{\tau}_{\bbb}\>$, where $\bbb = \{\<B_{\xi}^0,B_{\xi}^1\>
: {\xi}<\lambda\}$ is a separating and independent
family of $2$-partitions of ${\kappa}$. Moreover, $X$
is $\mu$-dense in ${D(2)}^{\lambda}$ iff $\bbb$ is $\mu$-independent.
\end{obs}

The spaces obtained from our main theorem
\ref{tm:main} will all be of the above form, with $\lambda =
2^\kappa$. The following fact will be instrumental in finding
appropriate families of dense sets $\dcal$ to be used to
produce $\dcal$-forced spaces.

\begin{fact}\label{tm:ld}
For each infinite cardinal ${\kappa}$, there is a
 family
$$\bbb=\{\<B_{\xi}^i:i<{\kappa}\>:{\xi}<2^{\kappa}\}$$
of $2^\kappa$ many ${\kappa}$-partitions of ${\kappa}$ that is
${\kappa}$-independent.
\end{fact}

Indeed, this fact is just a reformulation of the statement that
the space $\dis({\kappa})^{2^{\kappa}}$, the $2^\kappa$th power of
the discrete space on $\kappa$, contains a ${\kappa}$-dense subset
of size ${\kappa}$. This, in turn, follows immediately from the
Hewitt-Marczewski-Pondiczery theorem, see e. g. \cite[theorem
2.3.15]{En}.

\begin{mtheorem}\label{tm:main}
Assume that ${\kappa}$ is an infinite cardinal and we are given
$\bbb=
\bigl\{\<B_{\xi}^0,B_{\xi}^1\>
: {\xi}<2^{\kappa}
\bigr\}$, a ${\kappa}$-independent family of  $2$-partitions of
${\kappa}$,  moreover a non-empty family $\dcal$ of
${\kappa}$-dense subsets of the space $X_{\bbb}$. Then there is
another, always separating, ${\kappa}$-independent family
$\cbb=\{\<C_{\xi}^0,C_{\xi}^1\>:{\xi}< 2^{\kappa}\}$ of
$2$-partitions
 of ${\kappa}$  that satisfies the following five conditions:
\begin{enumerate}[(1)]
\item \label{dense}
every $D\in \dcal$ is also ${\kappa}$-dense in $X_{\cbb}$ (and so
$\operatorname{\Delta}(X_\cbb)={\kappa}$),
\item \label{forced} $X_{\cbb}$ is $\dcal$-forced,
\item \label{small}  $\nwd(X_{\cbb})={\kappa}$, i.e.
$\br {\kappa};<{\kappa};\subs \nwds(X_{\cbb})$,
\item \label{nowhere} $X_{\cbb}$ is  NODEC.
\newcounter{x}
\setcounter{x}{\value{enumi}}
 \end{enumerate}
Moreover, if $J\subs 2^{\kappa}$ is given with
$|2^{\kappa}\setm J|=2^{\kappa}$ then we can
assume that
\begin{enumerate}[(1)]
\addtocounter{enumi}{\value{x}}
\item \label{bj} $\cbb\restr J=\bbb \restr J$.
  \end{enumerate}
\end{mtheorem}

\begin{proof}

Assume that $J$ is given and let $I=2^{\kappa}\setm J$. We
partition $I$ into two disjoint pieces, $I=I_0\cup I'$, such that
$|I_0|={\kappa}^{<{\kappa}}$ and $|I'|=2^{\kappa}$.  Next we
partition $I_0$ into pairwise disjoint countable sets
$J_{A,{\alpha}}\in \br I_0;{\omega};$ for all $A\in \br
{\kappa};<{\kappa};$ and ${\alpha}\in {\kappa}\setm A$. If
${\xi}\in J_{A,{\alpha}}$ (for some $A\in \br {\kappa};<{\kappa};$
and ${\alpha}\in {\kappa}\setm A$) then we let
\begin{displaymath}
C^0_{\xi}=(B^0_{\xi}\cup A )\setm \{{\alpha}\},
\end{displaymath}
and
\begin{displaymath}
C^1_{\xi}=(B^1_{\xi}\setm A )\cup \{{\alpha}\}.
\end{displaymath}

Next, let us fix any enumeration $\{F_{\nu}:{\nu}<2^{\kappa}\}$
of $[\kappa]^{\kappa}$ and then by transfinite recursion on
${\nu}<2^{\kappa}$ define
\begin{itemize}
\item sets $K_{\nu}\subset  I'$ with $K_\nu = \emptyset$ or $|K_\nu| = \kappa$,
\item partitions $
\<C^0_{\sigma}, C^1_{\sigma}\>$ of ${\kappa}$ for all ${\sigma}\in
K_{\nu}$,
\item finite functions
${\eta}_{\nu}\in \fn(2^{\kappa},2)$,
\end{itemize}
such that the inductive hypothesis
\begin{equation}\tag{$\phi_{\nu}$}
\forall\varepsilon \in \fn(2^{\kappa},2)\ 
\forall D\in \dcal\ 
|D\cap \bbb_{{\nu}}[\varepsilon]|={\kappa}
\end{equation}
holds,
where
\begin{displaymath}
\bbb_{\nu}=\Bigl\{\<C^0_{\sigma},C^1_{\sigma}\>:{\sigma}\in
I_\nu\Bigr\}
\cup \Bigl\{\<B^0_{\sigma},B^1_{\sigma}\>:{\sigma}\in
2^{\kappa}\setm I_\nu\Bigr\}
\end{displaymath}
with
 $$I_\nu = I_0\cup\bigcup_{{\zeta}<{\nu}}K_{\zeta}.$$

Note that $(\phi_{\nu})$ simply says that every set $D \in \dcal$
is $\kappa$-dense in the space $X_{\bbb_{\nu}}$. We shall then
conclude that $\cbb=\bbb_{2^{\kappa}}$ is as required.

Let us observe first that $(\phi_0)$ holds because, by assumption,
we have $|\bbb[{\varepsilon}]\cap D|={\kappa}$ for all $D\in
\dcal$ and ${\varepsilon}\in \fn(2^{\kappa},2)$, moreover
\begin{displaymath}
\bigl|\bbb[{\varepsilon}]\bigtriangleup\bbb_0[{\varepsilon}]\bigr|<{\kappa}.
\end{displaymath}

Clearly, if ${\nu}$ is a limit ordinal and $(\phi_{\zeta})$ holds for
each ${\zeta}<{\nu}$ then $(\phi_{\nu})$ also holds. So the induction
hypothesis is preserved in limit steps.

Now consider the successor steps. Assume that $(\phi_{\nu})$
holds. We distinguish two cases:
\begin{case}
$F_{\nu}$ contains a \pieu{\dcal}{X_{\bbb_{\nu}}}, i.e.
$F_{\nu}\supset D\cap \bbb_{\nu}[{\eta}_{\nu}]$ for some
${\eta}_{\nu}\in \fn(2^{\kappa},2)$ and $D\in\dcal$.
\end{case}

This defines ${\eta}_{\nu}$ and then we set $K_{\nu}=\empt$. The
construction from here on will not change the partitions whose
indices occur in $\dom ({\eta}_{\nu})$, thus we shall have
$\bbb_{\nu}[{\eta}_{\nu}]=\bbb_{2^{\kappa}}[{\eta}_{\nu}]$ and so
at the end $F_{\nu}$ will include the
\pieu{\dcal}{X_{{\bbb_{2^{\kappa}}}}} $D \cap
\bbb_{2^{\kappa}}[{\eta}_{\nu}]$. Also, in this case we have
$\bbb_\nu = \bbb_{\nu+1}$, hence $(\phi_{\nu+1})$ trivially
remains valid.

\begin{case}$F_{\nu}$ does not include a \pieu{\dcal}{X_{\bbb_{\nu}}},
i.e. $(D\cap \bbb_{\nu}[{\varepsilon}])\setm F_{\nu}\ne \empt $
for all ${\varepsilon}\in \fn(2^{\kappa},2)$ and $D\in\dcal$.
\end{case}

In this case we
choose and fix any set
\begin{equation}\notag
K_{\nu}\subs I'\setm \bigl( \cup\{\dom {\eta}_{\zeta}:{\zeta}<
{\nu}\}\cup\cup
 \{K_{\zeta}:{\zeta}<{\nu}
\}
\bigr)
\end{equation}
of size ${\kappa}$ and let $K_\nu =
\{{\gamma}_{{\nu},i}:i<{\kappa}\}$ be a 1-1 enumeration of
$K_{\nu}$. We also set ${\eta}_{\nu}=\empt$. We want to modify the
partitions with indices in $K_\nu$ so as to make the set $F_\nu$
closed discrete in $X_{\bbb_{\nu+1}}$ and hence in
$X_{\bbb_{2^{\kappa}}}$ as well. To do this, we set for all
$i<{\kappa}$
\begin{displaymath}
C^0_{{\gamma}_{{\nu},i}}=(B^0_{{\gamma}_{{\nu},i}}\setm F_{\nu})\cup\{i\},
\end{displaymath}
and
\begin{displaymath}
C^1_{{\gamma}_{{\nu},i}}=(B^1_{{\gamma}_{{\nu},i}}\cup F_{\nu})\setm \{i\}.
\end{displaymath}
Then for each $i \in \kappa$ we have $i \in
C^0_{{\gamma}_{{\nu},i}}$ and $$F_\nu \cap
C^0_{{\gamma}_{{\nu},i}} \subset \{i\},$$ consequently $F_{\nu}$
is closed discrete in $X_{\bbb_{{\nu}+1}}$, hence $F_{\nu}$ will
be closed discrete in $X_{\bbb_{2^{\kappa}}}$.

We still have to show that 
$(\phi_{{\nu}+1})$ holds in this case, too. Assume, indirectly,
that for some $D\in\dcal$ and ${\varepsilon}\in \fn(2^{\kappa},2)$
we have
\begin{equation}\notag
\bigl|D\cap \bbb_{{\nu}+1}[{\varepsilon}]\bigr|<{\kappa}.
\end{equation}
Then we can clearly find ${\xi}\in I_0\setm \dom {\varepsilon}$
with
$$(D\cap \bbb_{{\nu}+1}[{\varepsilon}])\cup \dom(\varepsilon)\subs C^0_{\xi},$$
and so for ${\varepsilon}^*={\varepsilon}\cup \{\<{\xi},1\>\}$ we
even have
\begin{equation}\notag
D\cap \bbb_{{\nu}+1}[{\varepsilon}^*]=\empt.
\end{equation}
On the other hand, our choices clearly imply that
$$\bbb_{{\nu}+1}[{\varepsilon}^*] \supset
\bbb_{{\nu}}[{\varepsilon}^*]\setm F_{\nu},$$ consequently
$$D\cap
\bbb_{{\nu}+1}[{\varepsilon}^*] \supset (D\cap
\bbb_{{\nu}}[{\varepsilon}^*])\setm F_{\nu} \neq \emptyset,$$
a contradiction. This shows that $(\phi_{{\nu}+1})$ is indeed
valid, and the transfinite construction of $\cbb =
\bbb_{2^\kappa}$ is thus completed. We show next that $\cbb$
satisfies all the requirements of our main theorem.

$\cbb$ is \separ\ because e. g.
 for any ${\xi}\in J_{\{{\alpha}\},{\beta}}$ the partition
$\<C^0_{\xi}, C^1_{\xi}\>$ separates ${\alpha}$ and ${\beta}$.

That $\cbb$ is $\kappa$-independent and that (\ref{dense}) holds
(i. e. each $D\in \dcal$ is ${\kappa}$-dense in $X_{\cbb}$) both
follow from $(\phi_{2^{\kappa}})$.

If $A\in \br {\kappa};<{\kappa};$ and ${\alpha}\in {\kappa}\setm
A$, then  for any ${\xi}\in J_{A,{\alpha}}$ we have $A\subs
C^0_{\xi}$ and ${\alpha}\in C^1_{\xi}$, hence ${\alpha}\notin
\overline A^{X_{\cbb}}$. Thus every member of $ \br
{\kappa};<{\kappa};$ is closed and hence closed discrete in
$X_{\cbb}$, and so (\ref{small}) is satisfied.

Assume next that $F\in \ncal( X_{\cbb})$, we want to show that $F$
is closed discrete . By (\ref{small}) we may assume that
$|F|={\kappa}$ and so can find ${\nu}<2^{\kappa}$ with
$F=F_{\nu}$. Suppose that  at step ${\nu}$ of the recursion we
were in case 1; then we had $F\supset D\cap
\bbb_{\nu}[{\eta}_{\nu}]$ for some $D\in\dcal$. But
$\bbb_{\nu}[{\eta}_{\nu}]=\bbb_{2^{\kappa}}[{\eta}_{\nu}]=
\cbb[{\eta}_{\nu}]$, so $F$ would be dense in
$\cbb[{\eta}_{\nu}]$. This contradiction shows that, at step
$\nu$, we must have been in case 2. However, in this  case we know
that $F = F_{\nu}$ was made to be closed discrete in
$X_{\bbb_{\nu+1}}$ and consequently in $X_{\cbb}$ as well.
So $X_{\cbb}$ is NODEC, i.e. (\ref{nowhere}) holds.

It remains to check that $X_{\cbb}$ is $\dcal$-forced, i. e. that
(\ref{forced}) holds. By \ref{f:pd} it suffices to show that any
somewhere dense subset $E$ of $X_\cbb$ includes a
\pieu{\dcal}{X_{\cbb}}.
By (\ref{small}) we must have $|E| = \kappa$ and hence we can pick
${\nu}<2^{\kappa}$ such that $F_{\nu}= E$.
Then at step ${\nu}$ of the recursion we could not be in case 2,
since otherwise
 $F_{\nu} = E$ would have been made closed discrete in
$X_{\bbb_{{\nu}+1}}$ and so in $X_{\cbb}$ as well.
Hence at step ${\nu}$ of the recursion we were in case 1,
consequently ${\eta}_{\nu}\in \fn(2^{\kappa},2)$   and $D\in\dcal$
could be found such that $E = F_{\nu}\supset D\cap
\bbb_{\nu}[{\eta}_{\nu}]$.
However, by the construction, we have
$\cbb[{\eta}_{\nu}]=\bbb_{\nu}[{\eta}_{\nu}]$, and therefore $E$ actually
includes the \pieu{\dcal}{X_{\cbb }} $D\cap \cbb[{\eta}_{\nu}]$.

Finally, (\ref{bj}) trivially holds by the construction.

\end{proof}

\section{Applications to resolvability}

In this and the following section we shall present a large number
of consequences of our main theorem \ref{tm:main}. The key to most
of these will be given by a judicious choice of a family $\dcal$
of $\kappa$-dense sets in a space $X_{\bbb}$, where  $\bbb =
\bigl\{\<B_{\xi}^0,B_{\xi}^1\> : {\xi}<2^{\kappa} \bigr\}$ is a
${\kappa}$-independent family of $2$-partitions of some cardinal
${\kappa}$. In our first application, however, this choice is
trivial, that is we have $\dcal = \{\kappa\}$.

In \cite{ASTTW}, the following results were proven:
\begin{enumerate}[(1)]
\item  $\dis(2)^{\mathfrak c}$ does not have
a  dense countable  maximal subspace,
\item $\dis(2)^{\mathfrak c}$  has
a dense countable  irresolvable subspace,
\item  it is consistent that
$\dis(2)^{\mathfrak c}$ has a dense countable  submaximal
subspace,
  \end{enumerate}
and then the following natural problem was raised (\cite[Question
4.4]{ASTTW}): {\em Is it provable in ZFC that the Cantor cube
$\dis(2)^{\mathfrak c}$ or the Tychonov cube $[0,1]^{\mathfrak c}$
 has a dense countable submaximal subspace?}
Our next result gives an affirmative answer to this problem.

\begin{theorem}\label{tm:inte}
For each infinite cardinal ${\kappa}$ the Cantor cube
$\dis(2)^{2^{\kappa}}$ contains a dense submaximal  subspace $X$
 with  $|X|=\Delta(X)={\kappa}$.
\end{theorem}

\begin{proof}\prlabel{tm:inte}
Let us start by fixing any ${\kappa}$-independent family of
$2$-partitions $\bbb=\{\langle B^0_{\xi},B^1_{\xi}\rangle :
{\xi}<2^{\kappa}\}$ of ${\kappa}$, and let
$\dcal=\bigl\{{\kappa}\bigr\}$. Applying theorem \ref{tm:main}
with $\bbb$ and $\dcal$
 we obtain a family of 2-partitions $\cbb$ of $\kappa$ that satisfies
\ref{tm:main} (1)--(4). The space $X_{\cbb}$ is as required.
Indeed, $\Delta(X_{\cbb})={\kappa}$ because of
\ref{tm:main}(\ref{dense}), $X_{\cbb}$ is NODEC by
\ref{tm:main}(\ref{nowhere}), while it is OHI by lemma
\ref{i:irres}. But then it is submaximal. Finally, by observation
\ref{lm:embed},  $X_{\cbb}$ embeds  into $\dis(2)^{2^{\kappa}}$ as
a dense subspace.
\end{proof}

That theorem \ref{tm:inte} fully answers  \cite[Question
4.4]{ASTTW} follows from the following fact \ref{f:homeo} that may be
already known, although we have not found it in the literature.

 \begin{fact}\label{f:homeo}
 Any countable dense subspace of $\dis(2)^{\mathfrak c}$ is homeomorphic to a {\em
 dense} subspace of $[0,1]^{\mathfrak c}$.
 \end{fact}

This fact, in turn, immediately follows from the next proposition.
In it, as usual, we denote by $\irrac$ the space of the
irrationals.

 \begin{proposition}\label{p:homeo}
 Assume that $\kappa$ is an infinite cardinal, $S\subs \dis(2)^{\kappa}$ is dense, moreover there is a partition
$\{I_{\nu}:{\nu}<{\kappa}\}$ of ${\kappa}$ into countably infinite
sets such that for each ${\nu}<{\kappa}$ the set $2^{I_{\nu}}\setm
(S\restriction I_{\nu})$ is dense (in other words: $S\restriction
I_{\nu}$ is co-dense) in $2^{I_{\nu}}$. (The last condition is
trivially satisfied if the cardinality of $S$ is less than continuum.) Then $S$ is homeomorphic
to a dense subspace of the irrational cube $\irrac^{\kappa}$ and
hence of the Tychonov cube $[0,1]^\kappa$.
 \end{proposition}

 \begin{proof}
 For each ${\nu}<{\kappa}$ we may select a countable dense subset of $D_{\nu}\subs 2^{I_{\nu}}\setm
(S\restriction I_{\nu})$. The space $2^{I_{\nu}}\setm D_{\nu}$ is
known to be homeomorphic to $\irrac$ for all $\nu < \kappa$. Also,
for each $\nu < \kappa$ we have $S\restriction I_{\nu}\subs
2^{I_{\nu}}\setm D_{\nu}$ and therefore $S$ is naturally
homeomorphic to a dense subspace of the product space
$$\prod\{2^{I_{\nu}}\setm D_{\nu} : \nu < \kappa\}.$$ This
product, however, is homeomorphic to the cube $\irrac^{\kappa}$.
 \end{proof}

Let us remark that, as far as we know, the first ZFC example of a
countable regular, hence $0$-dimensional, submaximal space was
constructed by E. van Douwen in \cite{vD2}, by using an approach
that is very different from and much more involved than ours.
Also, it is not clear if his example embeds {\em densely} into the
Cantor or Tychonov cube of weight $\mathfrak c$.

After proving in \cite[Corollary 8.5]{AC} that every separable
submaximal topological group is countable,  Arhangel'skii and
Collins raised the following question \cite[Problem 8.6]{AC}: {\em
Is there a crowded uncountable separable Hausdorff (or even
Tychonov) submaximal space?} As it turns out, starting from any
zero-dimensional countable submaximal space (e. g. the one
obtained from the previous theorem or van Douwen's example from
\cite{vD2}) an affirmative answer can be given to this question,
at least in the $T_2$ case. The regular or Tychonov cases of the
problem, however, remain open.

\begin{theorem}\label{tm:subc}
There is a  crowded, separable, submaximal $T_2$ space $Y$  of
cardinality $\mathfrak c$. 
\end{theorem}

\begin{proof}
Let  $\tau$ be any fixed crowded, submaximal, $0$-dimensional, and
$T_2$ topology on $\omega$. Since $\tau$ is not compact we can
easily find $\{U_{\sigma}:{\sigma}\in 2^{<{\omega}}\}$, an
infinite partition of ${\omega}$ into  nonempty $\tau$-clopen sets
indexed by all finite 0-1 sequences $\sigma$.

The underlying set of $Y$ will be ${\omega}\cup{}^{\omega}2$ and
we let $X = \langle\omega,\tau\rangle$ be an open subspace of $Y$.
Next, a basic neighbourhood of a point $f\in {}^{\omega}2$ will be
of the form
\begin{displaymath}
\{f\}\cup\bigcup\{D_{f\restr n}:n\ge m\},
  \end{displaymath}
where $m\in {\omega}$ and $D_{f\restr n}$ is a dense (hence, as
$X$ is submaximal, open) subset of $U_{f\restr n}$ for $m\le
n<{\omega}$. It is easy to see that $Y$ is $T_2$, and $Y$ is
separable because $\omega$ is dense in it.

Now, assume that $D\subs Y$ is dense. Then $D\cap X$ is also dense
hence open in $X$, and similarly $D\cap U_{{\sigma}}$ is dense
open in $U_{\sigma}$ for each
 ${\sigma}\in 2^{<{\omega}}$. So for each $f\in D$  the set
$\{f\}\cup \bigcup\{D\cap U_{f\restr n}:n\ge 0\}\subs D$ is a
basic neighbourhood of $f$, showing that $D$ is open in Y.
\end{proof}

In 1967 Ceder and Pearson, \cite{CP}, raised the question whether
an ${\omega}$-resolvable space is necessarily  maximally
resolvable?
 El'kin,  \cite{El}, constructed a $T_1$ counterexample to this question, and then
Malykhin, \cite{Ma2}, produced  a crowded \hres\ $T_1$ space (that
is clearly ${\omega}$-resolvable) which is not maximally
resolvable.
 Eckertson, \cite{E}, and later Hu,
\cite{Hu}, gave  Tychonov counterexamples but not in ZFC:
Eckertson's construction used a  measurable cardinal, while Hu
applied the assumption $2^{\omega}=2^{\oo}$. Whether one could
find a Tychonov counterexample to the Ceder-Pearson problem in ZFC
was repeatedly asked as recently as in \cite{Co} and \cite{CoHu}.

Our next theorem gives a whole class of $0$-dimensional $T_2$
(hence Tychonov) counterexamples to the Ceder-Pearson problem in
ZFC. Quite naturally, they involve applications of our main
theorem \ref{tm:main} where the family of dense sets $\dcal$ forms
a partition of the underlying set.

Recall that any application of theorem \ref{tm:main} yields a
dense NODEC subspace $X$ of some Cantor cube $\dis(2)^{2^\kappa}$
with the extra properties $$|X| = \nwd(X) =
\operatorname{\Delta}(X) = \kappa.$$ From now on, we shall
call any space having all these properties a $\ccal(\kappa)$-space.
Of course, any $\ccal(\kappa)$-space is zero-dimensional $T_2$ and
therefore Tychonov. Finally, with the intention to use lemma
\ref{i:space}, we recall that any $\ccal(\kappa)$-space $X$, being dense in
a Cantor cube, is
CCC, i. e. satisfies $\widehat{c}(X) = \omega_1.$

\begin{theorem}\label{tm:mplu}
For any two infinite cardinals ${\mu}< {\kappa}$ there is a
$\ccal(\kappa)$-space $X$ that is the disjoint union of ${\mu}$
dense submaximal subspaces but is not ${\mu}^+$-\eres. (Of
course, $X$ is then $\mu$-{\res}  but not $\mu^+$-\res, 
hence not maximally resolvable.)
\end{theorem}

\begin{proof}
Using  fact \ref{tm:ld} we can easily find a ${\mu}$-partition
$\<D_i:i<{\mu}\>$ and  a  family of $2$-partitions
$\bbb=\bigl\{\<B_{\xi}^0,B_{\xi}^1\>:{\xi}<2^{\kappa}\bigr\}$ of
${\kappa}$ such that for each $i<{\mu}$ and ${\varepsilon}\in \fn
(2^{\kappa},2)$ we have
\begin{displaymath}
\bigl|\ D_i\cap \mathbb B[{\varepsilon}]\ 
\bigr|={\kappa}.
\end{displaymath}
We may then apply theorem \ref{tm:main} to this $\bbb$ and the
family $\dcal = \{D_i:i<{\mu}\}$
 to get a collection $\cbb$ of $2$-partitions  of ${\kappa}$ satisfying
\ref{tm:main}(1)-(4). We claim that the space $X_\cbb$ is as
required.

Firstly, as the members of $\dcal$  partition  ${\kappa}$ and
$X_{\cbb}$ is NODEC, lemma \ref{i:irres} implies that each $D_i
\in \dcal$ is a submaximal dense subspace of $X_{\cbb}$.

Secondly, since $X_{\cbb}$ is CCC and $|\dcal|={\mu} \ge \omega$,
lemma \ref{i:space}  implies that $X_{\cbb}$ is not
${\mu}^+$-\eres.

\end{proof}

Theorem \ref{tm:mplu} talks about infinite cardinals, and with
good reason; it has been long known that for any finite $n$ there
are say countable zero-dimensional spaces that are $n$-{\res}  but
not $(n+1)$-{\res}. In connection with this, Eckertson asked in
\cite[Question 4.5]{E} the following question: {\em Does there exist
for each infinite cardinal ${\kappa}$ and for each natural number
$n\ge 1$ a Tychonov space $X$ with
$|X|=\operatorname{\Delta}(X)={\kappa}$ such that $X$ is
$n$-resolvable but $X$ contains no $(n+1)$-resolvable subspaces?}
Li Feng, \cite{F}, gave a positive answer to this question  and
the following corollary of \ref{tm:mplu} improves his result. Our
example is a $\ccal(\kappa)$-space that is the disjoint union of
$n$ dense submaximal subspaces.

\begin{corollary}\label{tm:nno}
For each cardinal ${\kappa}\ge {\omega}$ and each natural number
$n\ge 1$ there is a $\ccal(\kappa)$-space $Y$
 which is the disjoint union of $n$ dense submaximal subspaces.
Then $Y$, automatically, does not contain any $(n+1)$-resolvable
subspaces.
\end{corollary}

\begin{proof}
Consider the $\ccal(\kappa)$-space $X$ given by  theorem
\ref{tm:mplu} for any fixed pair of cardinals $\mu < \kappa$ and
then set
$Y=\bigcup\{D_i:i<n\}$. Here each subspace $D_i$ of $Y$ is
submaximal and therefore HI. Consequently, every subspace of $Y$
can be written as the union of at most $n$ HI subspaces.
By \cite[lemma 2]{I}, no such space can be $(n+1)$-resolvable,
hence $Y$ contains no $(n+1)$-resolvable subspaces.

\end{proof}

Another question that can be raised concerning theorem
\ref{tm:mplu} is whether it could be extended to apply to all
infinite cardinals instead of just the successors $\mu^+$. It is
actually known that the answer to this question is negative.

Indeed, Illanes, and later Bashkara Rao proved the following two
``compactness''-type results on ${\lambda}$-resolvability, for
cardinals ${\lambda}$ of countable cofinality.

\begin{qtheorem}[Illanes, \cite{I}]
If a topological space $X$ is $k$-resolvable for each $k<{\omega}$
then $X$ is ${\omega}$-resolvable.
\end{qtheorem}

\begin{qtheorem}[Bhaskara Rao, \cite{Ba}]
If ${\lambda}$ is a singular cardinal with
$\cf({\lambda})={\omega}$ and $X$ is any topological space  that
is
  $\mu$-resolvable for each $\mu <\lambda$ then
$X$ is ${\lambda}$-resolvable.
\end{qtheorem}

In contrast to these, our next result,  theorem \ref{tm:limit},
implies that no such compactness-phenomenon is valid for
uncountable regular limit (that is inaccessible) cardinals.
However, the following intriguing problem remains open.
\begin{problem}
Assume that  ${\lambda}$ is a singular cardinal with
$\cf({\lambda})>{\omega}$ and $X$ is a topological space that is
${\mu}$-resolvable for all ${\mu}<{\lambda}$. Is it true then that
$X$ is also ${\lambda}$-resolvable?
\end{problem}

Theorem \ref{tm:limit} may be viewed as an extension of
\ref{tm:mplu} from successors to all uncountable regular
cardinals, providing counterexamples to the Ceder-Pearson problem
in further cases. However, the spaces obtained here are quite
different from the ones constructed in \ref{tm:mplu} because they
are { \em hereditarily resolvable}.

\begin{theorem}\label{tm:limit}
For any two cardinals $\kappa$ and $\lambda$ with
${\omega}<\cf({\lambda})={\lambda}\le {\kappa}$ there is a
$\ccal(\kappa)$-space that
is not ${\lambda}$-extraresolvable (and hence not
$\lambda$-resolvable) and still it is hereditarily
$\mu$-resolvable for all $\mu<{\lambda}$.
\end{theorem}

\begin{proof}
\prlabel{tm:limit}
Let us fix the sequence $\vec
{\lambda}=\<{\lambda}_{\zeta}:{\zeta}<{\lambda}\>$ by setting
${\lambda}_{\zeta}={\rho}$ for each ${\zeta}<{\lambda}$ if
${\lambda}={\rho}^+$ is a successor and by putting $\lambda_\zeta
= \omega_\zeta $ for ${\zeta}<{\lambda}$ if $\lambda$ is a limit
cardinal (note that $\lambda = \omega_\lambda$ in the latter
case).

Using  fact \ref{tm:ld} we can find two families of partitions
$$\mathbb
D=\bigl\{\<D_{\zeta}^i:i<{\lambda}_{\zeta}\>:{\zeta}<{\lambda}\bigr\}\\\ 
\mbox{and}\\\ 
 \mathbb
B=\bigl\{\<B_{\xi}^0,B_{\xi}^1:\>:{\xi}<2^{\kappa}\bigr\}$$
 of
${\kappa}$
such that $\dbb\cup \bbb$ is ${\kappa}$-independent, i. e. $\bigl|
\mathbb D[{\eta}]\cap \mathbb B[{\varepsilon}] \bigr|={\kappa}$
whenever ${\eta}\in \finc(\vec {\lambda})$ and ${\varepsilon}\in
\fn (2^{\kappa},2)$.
Then
\begin{displaymath}
\dcal=\{\dbb[{\eta}]:{\eta}\in \finc(\vec {\lambda})\}
\end{displaymath}
is a family of $\kappa$-dense sets in the space $X_{\bbb}$, hence
we can apply theorem \ref{tm:main} with $\bbb$ and $\dcal$
 to get a family $\cbb$ of $2$-partitions  of ${\kappa}$ satisfying
\ref{tm:main}(1)--(4). We shall show that the
$\ccal(\kappa)$-space $X_{\cbb}$ is as required.

\begin{claim}\label{cl:lambda}
For every family $\ecal\in \br\dcal;{\lambda};$ there is $\fcal\in
\br \ecal;{\lambda};$ such that $D\cap D' \in \dcal$ (and hence is
dense in $X_{\cbb}$) whenever $\{D,D'\}\in \br \fcal;2;$.
\end{claim}

\begin{proof}
\prlabel{cl:lambda}
We can write $\ecal =
\{\dbb[{\eta}_{\gamma}]:{\gamma}<{\lambda}\}$.
Since ${\lambda}=\cf({\lambda})>{\omega}$  we can find $K\in\br
{\lambda};{\lambda};$ such that
$\{\dom({\eta}_{\gamma}):{\gamma}\in K\}$ forms a $\Delta$-system
with  kernel $K^*$. Then $\prod_{i\in K^*}{\lambda}_i <{\lambda}$,
hence, as $\lambda$ is regular, there are a set $I\in \br
K;{\lambda};$ and a fixed finite function $\eta\in \prod_{i\in
K^*}{\lambda}_i \subset \finc(\vec {\lambda})$ such that
${\eta}_{\gamma}\restriction K^* ={\eta}$ for each ${\gamma}\in
I$.

But then $\fcal = \{\dbb[{\eta}_{\gamma}] : \gamma \in I\}$ is as
required: for any $\{{\gamma},{\delta}\}\in \br I;2;$ we have
${\eta}_{\gamma}\cup {\eta}_{\delta}\in \finc(\vec {\lambda})$,
consequently
$$\dbb[{\eta}_{\gamma}]\cap\dbb[{\eta}_{\delta}]=
\dbb[{\eta}_{\gamma}\cup {\eta}_{{\delta}}] \in \dcal.$$
\end{proof}

Now, since $\celh(X_{\cbb})=\oo$  and the above claim  holds we
can apply lemma \ref{i:space} to conclude that $X_{\cbb}$ is not
${\lambda}$-extraresolvable.

Let us now fix $\mu<{\lambda}$. We first show that  every
$\dbb[{\eta}] \in \dcal $ is $\mu$-resolvable. Indeed, choose
${\zeta}\in {\lambda}\setm\dom {\eta}$ with ${\lambda}_{\zeta}\ge
\mu$. Clearly, then the family $\{\dbb[{\eta}\cup
\{\langle\zeta,{\gamma}\rangle\}]:{\gamma}<{\lambda}_{\zeta}\}$
forms a partition of $\dbb[{\eta}]$ into ${\lambda}_{\zeta}\ge
\mu$ many dense subsets.

 Since
$X_{\cbb}$ is NODEC and $\dcal$-forced, any crowded subspace $S$
of $X_{\cbb}$ is somewhere dense.
Consequently, lemma 2.10 implies
that $X_{\cbb}$ is hereditarily $\mu$-resolvable.
\end{proof}

\begin{remark}
It is well-known that any dense subspace of the Cantor cube
$\dis(2)^{\lambda}$ has weight (even $\pi$-weight) equal to
$\lambda$. Consequently, any $\ccal(\kappa)$-space (that is, by
definition, of cardinality $\kappa$) has maximum possible weight,
that is $2^\kappa.$
Now, ZFC counterexamples to the Ceder-Pearson problem are
naturally expected to have this property.
Indeed, for instance the forcing axiom BACH (see e.g. \cite{Ta})
implies that every topological space $X$ with
$|X|=\operatorname{\Delta}(X)=\oo$ and $\piweight(X) <
2^{{\omega}_1}$ is $\oo$-resolvable.
Consequently, under BACH, any
$\omega$-resolvable space $X$ satisfying $|X| = \omega_1$ and
$\piweight(X) < 2^{{\omega}_1}$ is maximally resolvable.
\end{remark}

By \cite[Lemma 4]{I},  any topological space that is not
${\omega}$-resolvable contains a HI somewhere dense subspace.
 Theorem \ref{tm:limit} shows that this badly
fails if $\omega$ is replaced by an uncountable cardinal.

 Again by \cite[Lemma 4]{I}, if a space $X$ can be partitioned into finitely many dense HI
subspaces, then the number of pieces is uniquely determined. It
follows from our next result,  theorem \ref{tm:max} below, that
this is not the case for infinite partitions. In fact, for every
infinite cardinal ${\kappa}$ there is a $\ccal(\kappa)$-space that
can be simultaneously partitioned into ${\lambda}$ many dense
submaximal (and so HI) subspaces for all infinite ${\lambda}\le
{\kappa}$.

Theorem \ref{tm:max} also gives an affirmative answer to the
following question of Eckertson, raised in \cite[3.4 and 3.6]{E}:
{\em Does there exist, for each cardinal $\mu$, a
${\mu}^+$-resolvable space that can be partitioned into
${\mu}$-many dense HI subspaces? }

 The proof of theorem \ref{tm:max} will require an even more delicate choice of the family of dense sets
$\dcal$ than the one we used in the proof of \ref{tm:limit}.

\begin{theorem}\label{tm:max}
For each infinite cardinal ${\kappa}$ there is a
$\ccal(\kappa)$-space that can be simultaneously partitioned into
countably many dense hereditarily ${\kappa}$-resolvable subspaces
and also into $\mu$ many dense submaximal (and therefore HI)
subspaces for all infinite $ {\mu}\le {\kappa}$.
\end{theorem}

\begin{proof}\prlabel{tm:max}
Let us start by setting ${\lambda}_0={\omega}$,
${\lambda}_1={\kappa}$, and $\vec {\lambda}=\<{\lambda}_i:i<2\>$,
moreover
  $\vec {\kappa}=\<{\kappa}_n:n<{\omega}\>$,
where ${\kappa}_0={\omega}$ and ${\kappa}_n={\kappa}$ for $1\le
n<{\omega}$.

By fact \ref{tm:ld} there are three families of partitions of
${\kappa}$, say
$$\bbb=\{\<B_{\zeta}^i:i<2\>:{\zeta}<2^{\kappa}\},$$
$$\ebb=\{\<E_n^j:j<{\kappa}_n\>:n<{\omega}\},$$ and
$$\fbb=\{\<F_\ell^k:k<{\lambda}_\ell\>:\ell<2\},$$  such that
$\bbb\cup\ebb\cup\fbb$ is ${\kappa}$-independent, i.e.
for each ${\varepsilon}\in \fn(2^{\kappa},2)$, ${\eta}\in\ 
\finc(\vec {\kappa})$, and ${\rho}\in \finc (\vec {\lambda})$ we
have
  \begin{equation}\tag{\dag}\label{bef}
   \big |\bbb[{\varepsilon}]\cap \ebb[{\eta}]\cap \fbb[{\rho}]\big|={\kappa}.
  \end{equation}

Of course, (\ref{bef}) implies that all sets of the form
$\ebb[{\eta}]\cap \fbb[{\rho}]$ are $\kappa$-dense in $X_{\bbb}$,
however the family $\dcal$ of $\kappa$-dense sets that we need
will be defined in a more complicated way.

To start with, let us write
$\fcal_\ell=\{F^k_\ell:k<{\lambda}_\ell\}$ for $\ell<2$ and then
set
\begin{displaymath}
 \dcal_\ebb=\bigl\{\ebb[{\eta}]:{\eta}\in \finc(\vec {\kappa})\bigr\}
\end{displaymath}
and
\begin{displaymath}
 \dcal_\fbb = \fcal_0 \cup \fcal_1 =\{F_\ell^k:\ell<2,\\ k<{\lambda}_\ell\}.
\end{displaymath}
Next let
\begin{displaymath}\notag
\dcal_{\ebb,\fbb}=\{E\setm\cup \fcal: E\in \dcal_\ebb,\\  \fcal\in
\br \dcal_\fbb;<{\omega};\}
\end{displaymath}
and
\begin{displaymath}\notag
\dcal_{\fbb,\ebb}= \{F_\ell^k\setm \bigl((\cup \ecal)\cup(\cup
\fcal)\bigr): F_\ell^k\in \dcal_{\fbb},\\ \ecal\in \br
\dcal_\ebb;<{\omega};,\\  \fcal\in \br\fcal_{1-\ell};<{\omega};\}.
\end{displaymath}
Finally, we set
\begin{displaymath}
\dcal=\dcal_{\ebb,\fbb}\cup \dcal_{\fbb,\ebb}.
\end{displaymath}

Every  element of $\dcal$ contains some (in fact, infinitely many)
sets of the form $\ebb[{\eta}]\cap \fbb[{\rho}]$ and so is
 ${\kappa}$-dense in $X_{\bbb}$ by
(\ref{bef}).

Now we may apply theorem \ref{tm:main} with $\bbb$ and  $\dcal$ to
obtain a family of partitions $\cbb$ of $\kappa$ that satisfies
\ref{tm:main} (1) - (4). We shall show that $X_{\cbb}$ is as
required.

\begin{claim}\label{cl:nwd}
$E\cap F$ is nowhere dense in $X_{\cbb}$ whenever
$E\in \dcal_\ebb$ and $F\in \dcal_\fbb$.
\end{claim}

\begin{proof}\prlabel{cl:nwd}
According to  
\ref{i:nwd} it suffices to show that $D\setm (E\cap F)$ includes
an element of $\dcal$  whenever $D\in \dcal$.

Now, if $D=E'\setm \cup \fcal\in \dcal_{\ebb,\fbb}$ then
\begin{displaymath}
D\setm (E\cap F)\supset E'\setm (\cup(\fcal\cup\{F\}))\in \dcal
_{\ebb,\fbb}.
\end{displaymath}

If, on the other hand, $D=F_\ell^k\setm \bigl((\cup
\ecal)\cup(\cup \fcal)\bigr) \in \dcal_{\fbb,\ebb}$ then
\begin{displaymath}
D\setm (E\cap F)\supset
F_\ell^k\setm \bigl((\cup (\ecal\cup\{E\}))\cup(\cup \fcal)\bigr)
\in \dcal_{\fbb,\ebb}.
\end{displaymath}
\end{proof}

\begin{claim}\label{cl:nwd2}
$F\cap F'$ is nowhere dense in $X_{\cbb}$ for all
$\{F,F'\}\in \br \dcal_\fbb;2;$.
\end{claim}

\begin{proof}\prlabel{cl:nwd2}
Again, by  
\ref{i:nwd}, it is enough to show that $D\setm (F\cap F')$
includes an element of $\dcal$  for each $D\in \dcal$.

If $D=E\setm \cup \fcal\in \dcal_{\ebb,\fbb}$ then
\begin{displaymath}
D\setm (F\cap F')\supset E\setm (\cup(\fcal\cup\{F\}))\in \dcal
_{\ebb,\fbb}.
\end{displaymath}
If $D=F_\ell^k\setm \bigl((\cup \ecal)\cup(\cup \fcal)\bigr) \in
\dcal_{\fbb,\ebb}$ and $F\cap F'\ne \empt$ then  we can assume
that $F\in \fcal_{\ell}$ and $F'\in \fcal_{1-\ell}$. But then we
have
\begin{displaymath}
D\setm (F\cap F')\supset F_\ell^k\setm \bigl((\cup \ecal)\cup(\cup
(\fcal\cup\{F'\}))\bigr) \in \dcal_{\fbb,\ebb}.
\end{displaymath}
\end{proof}

\begin{claim}\label{cl:deres}
Every $D\in \dcal_{\ebb,\fbb}$ is ${\kappa}$-resolvable.
\end{claim}

\begin{proof}\prlabel{cl:deres}
Let $D=E\setm \cup\fcal$. Without loss of generality we can assume
that $E=\ebb[{\eta}]$ with $\dom {\eta}=n\in {\omega}\setm \{0\}$.
But then $D$ is the disjoint union of the $\kappa_n = \kappa$ many
dense sets
$$\{\ebb[{\eta} \cup \{\langle n,\zeta \rangle\}]\setm \cup \fcal :
{\zeta}<{\kappa}\}.$$
\end{proof}

\begin{claim}\label{cl:ekappa}
$E_0^i$ is hereditarily ${\kappa}$-resolvable for each $i<{\omega}
= \kappa_0 $.
\end{claim}

\begin{proof}
\prlabel{cl:ekappa}
Let us note first of all that for any $$D=F\setm
((\cup\ecal)\cup(\cup\fcal))\in  \dcal_{\fbb,\ebb}$$ we have
$E_0^i\cap D\subs E_0^i\cap F \in \ncal(X_{\cbb})$ by claim
\ref{cl:nwd}.

Now, let $S$ be any crowded subspace of $E_0^i$. Since $X_{\cbb}$
is NODEC and $\dcal$-forced, by lemma 2.10
there is a partial $(\dcal, X_{\cbb})$-mosaic $$M=\bigcup\{V\cap
D_V:V\in \vcal\} \subset S$$ that is dense in $S$. By our above
remark, we must have $D_V \in \dcal_{\ebb,\fbb}$ whenever $V \in
\vcal$, consequently $M$ and hence $S$ is ${\kappa}$-resolvable by
claim \ref{cl:deres} and fact \ref{i:res}.

\end{proof}

We have thus concluded that $\{E^i_0:i<{\omega}\}$  partitions
$X_{\cbb}$ into countably many hereditarily ${\kappa}$-resolvable
dense subspaces.

\begin{claim}\label{cl:fir}
$F_\ell^k \subset X_{\cbb}$ is submaximal for all $\ell<2$ and
$k<{\lambda}_\ell$.
\end{claim}

\begin{proof}
\prlabel{cl:fir}
Since $X_{\cbb}$ is NODEC, so is its dense subspace $F_{\ell}^k$,
hence it suffices to show that $F_\ell^k$ is OHI. By lemma
\ref{i:irres}, this will follow if we can show that for each
$D\in\dcal$ either $ F_\ell^k \cap D$ or $F_\ell^k\setm D$ is
nowhere dense in $X_{\cbb}$.

\newcases
\begin{case}
$D=E\setm \cup \fcal\in \dcal_{\ebb,\fbb}$.
\end{case}

Then $D\cap F_{\ell}^k\subs E\cap F_{\ell}^k \in \ncal(X_{\cbb})$
by claim \ref{cl:nwd}.

\begin{case}
$D=F'\setm ((\cup\ecal)\cup(\cup\fcal))\in  \dcal_{\fbb,\ebb}$.
\end{case}

If $F'\ne F_\ell^k$ then $F_{\ell}^k\cap D\subs F_{\ell}^k\cap F'
\in \ncal(X_{\cbb})$ by claim \ref{cl:nwd2}. Thus we may assume
that $F'=F_{\ell}^k$ and hence $F_{\ell}^k\notin \fcal$ because
$\fcal\subs \fcal_{1-\ell}$. But then
\begin{multline}\notag
F_{\ell}^k\setm D=
F_{\ell}^k\setm \big(F_{\ell}^k\setm
((\cup\ecal)\cup(\cup\fcal))\big)=\\
F_{\ell}^k\cap \big((\cup\ecal)\cup(\cup\fcal)\big)=
\cup_{E\in \ecal}(F_{\ell}^k\cap E)
\cup \cup_{F\in\fcal}(F\cap F_{\ell}^k),
\end{multline}
where each $F_{\ell}^k\cap E$ is nowhere dense by claim \ref{cl:nwd}
and each $F\cap F_{\ell}^k$ is nowhere dense by claim \ref{cl:nwd2}, i.e.
$F^k_\ell\setm D\in \nwds(X_{\cbb})$.
\end{proof}

Claim \ref{cl:fir} implies that $X_{\cbb}$ can be  partitioned
into $\mu$
many dense submaximal subspaces for both $\mu = \omega$ and $\mu =
\kappa.$
Since $\ccal(\kappa)$-spaces are CCC, it follows from theorem
\ref{tm:gap} below that this is also valid for all $\mu$ with
  ${\omega}<{\mu}<{\kappa}$.
\end{proof}

The following result is somewhat different from the others in that
it has no relevance to $\dcal$-forced spaces. Still we decided to
include it here not only because it makes the proof of theorem
\ref{tm:max} simpler but also because it seems to have independent
interest.

\begin{theorem}\label{tm:gap}
Let ${\omega}\le{\lambda}<{\mu}<{\kappa}$ be cardinals and $X$ be
a topological space with $\operatorname{c}(X)\le {\mu}$. If $X$
can be partitioned into both ${\lambda}$ many and ${\kappa}$ many
dense OHI subspaces then $X$ can also be partitioned  into ${\mu}$
many dense OHI subspaces.
\end{theorem}

\begin{proof}\prlabel{tm:gap}Let
$\<Y_{\sigma}:{\sigma}<{\lambda}\>$ and
$\<Z_{\zeta}:{\zeta}<{\kappa}\>$ be two partitions of $X$ into OHI
subspaces. For each ${\sigma}<{\lambda}$ let
\begin{multline}\notag
\ucal_{\sigma}=\{U\subs X:U \text{ is open and there is } \ 
{I}_{{\sigma},U}\in \br {\kappa};{\mu};\text{ such that }\\
Y_{\sigma}\cap\cup\{Z_{\zeta}:{\zeta}\in I_{{\sigma},U}\}\text{ is dense in }U\}.
\end{multline}
Since $\operatorname{c}(X)\le {\mu}$  there is
$\ucal_{\sigma}^*\in \br \ucal_{\sigma};\le{\mu};$ such that
$U_{\sigma}=\cup\ucal_{\sigma}^*$ is dense in
$\cup\ucal_{\sigma}$. Clearly, we also have $U_{\sigma}\in
\ucal_{\sigma}$. Next we set $V_{\sigma}=X\setm
\overline{U_{\sigma}}$ and
 $Q_{\sigma}=X\setm (U_{\sigma}\cup V_{\sigma})=\operatorname{Fr}(U_{\sigma})$.

Since $\lambda < \mu$ we can pick $I\in \br {\kappa};{\mu};$ with
\begin{displaymath}
\cup\{I_{{\sigma},U_{\sigma}}:{\sigma}<{\lambda}\} \subset I
\end{displaymath}
and then can choose  $J\in \br {\kappa}\setm I;{\lambda};$. Let
$Z=\bigcup\{Z_{\zeta}:{\zeta}\in I\cup J\}$.

For  ${\sigma}\in {\lambda}$ let $R_{\sigma}=Y_{\sigma}\cap
V_{\sigma}\cap Z$. Since $|I\cup J|={\mu}$, it follows from the
definition of $\ucal_\sigma$ and $V_{\sigma}=X\setm \overline{\cup
\ucal_{\sigma}}$ that
\begin{enumerate}
\item[($*$)]  $R_{\sigma}$ is nowhere dense in $X$ for each ${\sigma}<{\lambda}$.
\end{enumerate}

Let $P_{\sigma}=(Y_{\sigma}\cap U_{\sigma})\setm
\cup\{Z_{\zeta}:{\zeta}\in I_{{\sigma},U_{\sigma}}\}$ for
${\sigma}<{\lambda}$. Then $P_{\sigma}$ is also nowhere dense
because $\cup\{Z_{\zeta}:{\zeta}\in I_{{\sigma},U_{\sigma}}\}\cap
U_{\sigma}\cap Y_{\sigma}$ is dense in $U_{\sigma}$ and
$Y_{\sigma}$ is OHI.

Now let $\{{\sigma}_{\zeta}:{\zeta}\in J\}$ be an enumeration of
${\lambda}$ without repetition and  for each ${\zeta}\in J$ set
\begin{displaymath}
T_{{\zeta}}
=
(Z_{{\zeta}}\cap U_{{\sigma}_{\zeta}})\cup
\bigl((Y_{{\sigma}_{\zeta}}\cap V_{{\sigma}_{\zeta}})\setm Z\bigr)
.
\end{displaymath}

\begin{claim}
Each $T_{\zeta}$ is a dense OHI subspace of $X$.
\end{claim}

\begin{proof}
\prnolabel
$Z_{\zeta}$ is dense in $U_{{\sigma}_{\zeta}}$ and
$$(Y_{{\sigma}_{\zeta}}\cap V_{{\sigma}_{\zeta}})\setm Z=
(Y_{{\sigma}_{\zeta}}\cap V_{{\sigma}_{\zeta}})\setm
R_{{\sigma}_{\zeta}}$$ is dense in $V_{{\sigma}_{\zeta}}$ because
$Y_{{\sigma}_{\zeta}}$ is dense and
$R_{\sigma_{\zeta}}=Y_{{\sigma}_{\zeta}}\cap
V_{{\sigma}_{\zeta}}\cap Z$ is nowhere dense by $(*)$. Hence
$T_{\zeta}$ is dense. $T_{\zeta}$ is OHI because both $Z_{\zeta}$
and $Y_{{\sigma}_{\zeta}}$ are.
\end{proof}

\begin{claim}
The family $\{Z_{\xi}:{\xi}\in I\}\cup\{T_{\zeta}:{\zeta}\in J\}$
is disjoint.
\end{claim}

\begin{proof}
\prnolabel
Assume first that ${\xi}\in I$ and ${\zeta}\in J$. Then $\xi \ne
\zeta$ and hence

\begin{multline*}
T_{\zeta}\cap Z_{\xi}=\bigl((Z_{{\zeta}}\cap
U_{{\sigma}_{\zeta}})\cup \bigl((Y_{{\sigma}_{\zeta}}\cap
V_{{\sigma}_{\zeta}})\setm Z\bigr)\bigr)\cap Z_{\xi}\\
\subs (Z_{\zeta}\cap Z_{\xi})\cup (Z_{\xi}\setm Z)=\empt.$$
\end{multline*} Next if $\{{\zeta},{\xi}\}\in \br J;2;$, then

\begin{multline}\notag
T_{\zeta}\cap T_{\xi}=\\
\Bigl((Z_{{\zeta}}\cap U_{{\sigma}_{\zeta}})\cup
\bigl((Y_{{\sigma}_{\zeta}}\cap V_{{\sigma}_{\zeta}})\setm Z\bigr)\Bigr)\cap
\Bigr((Z_{{\xi}}\cap U_{{\sigma}_{{\xi}}})\cup
\bigl((Y_{{\sigma}_{{\xi}}}\cap V_{{\sigma}_{{\xi}}})\setm Z\bigr)\Bigr)\subs\\
(Z_{\zeta}\cap Z_{\xi})\cup (Z_{\zeta}\setm Z)\cup (Z_{\xi}\setm Z)\cup
(Y_{{\sigma}_{\zeta}}\cap Y_{{\sigma}_{\xi}})=\empt.
\end{multline}
\end{proof}

Thus we would be finished if we could prove that
$$\{Z_{\xi}:{\xi}\in I\}\cup\{T_{\zeta}:{\zeta}\in J\}$$ covers $X$.
However, we can only prove the following weaker statement.

\begin{claim}
$$X=\bigcup\{Z_{\xi}:{\xi}\in I\}\cup\bigcup\{T_{\zeta}:{\zeta}\in J\}
\cup\bigcup\{P_{\sigma}\cup Q_{\sigma}\cup
R_{\sigma}:{\sigma}<{\lambda}\}.$$
\end{claim}

\begin{proof}
\prnolabel
Let $x \in X$ be any point then there is a unique $\sigma <
\lambda$ with $x\in Y_{\sigma}$. If $x\notin U_{\sigma}\cup
V_{\sigma}$ then, by definition, $x\in Q_{\sigma}$.

So assume now that $x\in U_{\sigma}$. If $x\notin
\cup\{Z_{\zeta}:{\zeta}\in I_{{\sigma},U_{\sigma}}\}$ then $x\in
P_{\sigma}$. Otherwise $x\in Z_{\zeta}$ for some ${\zeta}\in
I_{{\sigma},U_{\sigma}}\subs I$.

Finally, assume that $x\in V_{\sigma}$ and let ${\zeta}\in J$ with
${\sigma}_{\zeta}={\sigma}$. Now, if $x\notin Z $ then $x\in
T_{\zeta}$ and if $x\in Z$ then $x\in R_{\sigma}$.
\end{proof}

The pairwise disjoint dense OHI subspaces $\{Z_{\xi}:{\xi}\in
I\}\cup \{T_{\zeta}:{\zeta}\in  J\}$ thus cover $X$ apart from the
nowhere dense sets $P_{\sigma}\cup Q_{\sigma}\cup R_{\sigma}$ for
$\sigma < \lambda.$ But then, using the obvious fact that the
union of a dense OHI subspace with any nowhere dense set is OHI,
the latter can be simply ``absorbed'' by the former, and thus a
partition of $X$ into $\mu$ many dense OHI subspaces can be
produced.
\end{proof}

\section{Applications to extraresolvability}

In \cite{CoHu} Comfort and Hu investigated the following question:
 {\em Are maximally resolvable spaces (strongly) extraresolvable?}
They presented several counterexamples, but the following 
question was left open (see \cite[Discussion
 1.4]{CoHu}):  
{\em Is there a maximally resolvable Tychonov space $X$
with $|X|=\nwd(X)$ such that $X$ is not extraresolvable?}
Using our main theorem \ref{tm:main} we can give an affirmative
answer to this question in ZFC. Recall that if $X$ is a
$\ccal(\kappa)$-space then $|X|=\nwd(X) = \kappa$.

\begin{theorem}\label{tm:extra}
For every infinite cardinal ${\kappa}$ there is a
$\ccal(\kappa)$-space
that  is  hereditarily ${\kappa}$-resolvable (and hence maximally
resolvable) but not \eres.
\end{theorem}

\begin{proof}

Let $\vec {\kappa}=\<{\kappa},{\kappa},\dots\>$ be the constant
${\kappa}$ sequence of length ${\omega}$. By  fact \ref{tm:ld}
there are a countable family $\mathbb
D=\bigl\{\<D_m^i:i<{\kappa}\>:m<{\omega}\bigr\}$ of
${\kappa}$-partitions of ${\kappa}$
 and a family $\bbb=\bigl\{\<B_{\xi}^0,B_{\xi}^1:\>:{\xi}<2^{\kappa}\bigr\}$
 of $2$-partitions
of ${\kappa}$ such that  $\bbb\cup\dbb$ is ${\kappa}$-independent,
that is
for each ${\eta}\in \finc(\vec {\kappa}) = \fn(\omega,\kappa)$ and
${\varepsilon}\in \fn (2^{\kappa},2)$ we have
\begin{displaymath}
\bigl|\ \mathbb D[{\eta}]\cap \mathbb B[{\varepsilon}]\ 
\bigr|={\kappa}.
\end{displaymath}

Now let
\begin{displaymath}
\dcal=\{\dbb[{\eta}]:{\eta}\in \finc(\vec {\kappa})\}
\end{displaymath}
and apply theorem \ref{tm:main} to $\bbb$ and $\dcal$ to get a
family $\cbb$ of $2$-partitions of ${\kappa}$  satisfying
\ref{tm:main} (1) - (4).

Since $|\dcal|={\kappa}$ and $\celh(X_{\cbb})=\oo$, it follows
from  lemma \ref{i:space} that $X_{\cbb}$ is not
${\kappa}^+$-\eres ( =  \eres).

Next, if $\dbb[{\eta}]\in \dcal$ then $\{\dbb[{\eta}\conc
\<{\zeta}\>]:{\zeta}<{\kappa}\}$  partitions  $\dbb[{\eta}]$ into
${\kappa}$ many dense sets, i.e. $\dbb[{\eta}]$ is
${\kappa}$-resolvable. Hence, by lemma 2.10,
$X_{\cbb}$ is hereditarily ${\kappa}$-resolvable.
\end{proof}

Our next two results are natural analogues of theorems
\ref{tm:mplu} and \ref{tm:limit} with ${\mu}$-resolvability
replaced by ${\mu}$-extraresolvability. Before formulating them,
however, we need a new piece of notation.

\begin{definition}\label{df:hint}
Given a family $\dbb=\bigl\{\<D^0_{\xi},D^1_{\xi}\>:{\xi}\in
{\rho}\bigr\}$ of $2$-partitions of a cardinal ${\kappa}$ we set
\begin{displaymath}
 \hint{\dbb}=\{D_{\zeta}^0\setm \cup_{{\xi}\in\Xi}D_{\xi}^0:
{\zeta}\in {\rho}\land
\Xi\in\br {\rho}\setm \{{\zeta}\};<{\omega};\}.
\end{displaymath}
\end{definition}

\begin{theorem}\label{tm:meres1}
For any infinite cardinals ${\kappa}\le {\lambda}\le 2^{\kappa}$
there is a ${\lambda}$-extraresolvable $\ccal(\kappa)$
 
-space $X$
that is not ${\lambda}^+$-\eres.  Moreover, every crowded subspace
of $X$ has a dense submaximal subspace.
\end{theorem}

\begin{proof}
\prtxtlabel{\ref{tm:meres1}}

By  fact \ref{tm:ld} there are families of $2$-partitions of
${\kappa}$, say $\mathbb D=\{\<D_{\zeta}^0,
D_{\zeta}^1\>:{\zeta}<{\lambda}\}$ and $\mathbb B=\{\<B_{\xi}^0,
B_{\xi}^1\>:{\xi}<2^{\kappa}\}$, such that $\bbb\cup\dbb$ is
${\kappa}$-independent, i. e.
$\bigl|\mathbb D[{\eta}]\cap \mathbb B[{\varepsilon}]
\bigr|={\kappa}$ for all ${\eta}\in \fn({\lambda},2)$ and
${\varepsilon}\in \fn (2^{\kappa},2)$.

Then
$\dcal=\hint{\dbb}$ is a family of $\kappa$-dense subsets of
$X_\bbb$,
hence we can apply the main theorem \ref{tm:main} to $\bbb$ and
$\dcal$ to obtain a family of partitions $\cbb$ satisfying
\ref{tm:main} (1) - (4).  We shall show that $X_{\cbb}$ is as
required.

\begin{claim}\label{cl:nwd3}
$D_{\zeta}^0\cap D_{\xi}^0\in \nwds(X_{\cbb})$ for each
pair $\{{\zeta},{\xi}\}\in \br {\lambda};2;$.
\end{claim}

\begin{proof}
\prnolabel
Write $Y=D_{\zeta}^0\cap D_{\xi}^0$ and  $D=D_{\nu}^0\setm
\cup_{{\eta}\in\Xi}D_{\eta}^0$ be an arbitrary member of $\dcal$.
 We can assume that ${\xi}\ne {\nu}$ and so
\begin{displaymath}
 D\setm Y=(D_{\nu}^0\setm \cup_{{\eta}\in\Xi}D_{\eta}^0)\setm
(D_{\zeta}^0\cap D_{\xi}^0)
\supset
D_{\nu}^0\setm \cup_{{\eta}\in\Xi\cup\{{\xi}\}}D_{\eta}^0\in\dcal,
  \end{displaymath}
showing that $D\setm Y$ is dense in $X_{\cbb}$. Hence, by  lemma
\ref{i:nwd}, $Y$ is nowhere dense in $X_{\cbb}$.
\end{proof}

Thus the family  $\{D_{\xi}^0:{\xi}\in {\lambda}\}$ witnesses that
$X_{\cbb}$ is ${\lambda}$-\eres. On the other hand, since
$|\dcal|={\lambda}$ and $\cel(X_{\cbb})={\omega}$, lemma
\ref{i:space} implies that $X_{\cbb}$ is not ${\lambda}^+$-\eres.

\begin{claim}\label{cl:sub}
Every $S\in \dcal$ is a submaximal subspace of $X_{\cbb}$.
\end{claim}

\begin{proof}
\prnolabel
Let $S=D_{\nu}^0\setm \cup_{{\eta}\in\Xi}D_{\eta}$, moreover
$D=D_{\mu}^0\setm \cup_{{\eta}\in\Psi}D_{\eta}^0$
 be an arbitrary member of $\dcal$.
If ${\nu}={\mu}$ then, by claim \ref{cl:nwd3},
\begin{displaymath}
S\setm D =(D_{\nu}^0\setm \cup_{{\eta}\in\Xi}D_{\eta}^0)\setm
(D_{\nu}^0\setm \cup_{{\eta}\in\Psi}D_{\eta}^0)\subs
\bigcup_{{\eta}\in \Psi}D_{\nu}^0\cap D_{\eta}^0\in
\nwds(X_{\cbb})
\end{displaymath}
and so $S\subs^* D$. If, on the other hand, ${\nu}\ne {\mu}$ then
we have
\begin{displaymath}
S\cap D= (D_{\nu}^0\setm \cup_{{\eta}\in\Xi}D_{\eta}^0)\cap
(D_{\nu}^0\setm \cup_{{\eta}\in\Psi}D_{\eta}^0)\subs D_{\nu}^0\cap
D_{\mu}^0\in \nwds(X_{\cbb})
\end{displaymath}
by claim \ref{cl:nwd3} again, consequently $S\cap D=^*\empt$. Thus
$S$ is OHI by lemma \ref{i:irres}, and since $X_{\cbb}$ is NODEC,
$S$ is even submaximal.
\end{proof}

Claim \ref{cl:sub} clearly implies that all $\dcal$-pieces and
hence all partial \mosaic{\dcal}s are submaximal subspaces of
$X_{\cbb}$. But $X_{\cbb}$ is $\dcal$-forced and NODEC, and
therefore, by lemma 1.10, every crowded subspace of $X_{\cbb}$
includes a partial \mosaic{\dcal} as a dense subspace.
\end{proof}

Let us remark that theorem \ref{tm:meres1} makes sense, and
remains valid, for $\lambda < \kappa$ as well. However, in this
case theorem \ref{tm:mplu} yields a stronger result. This is the
reason why we only formulated it for $\lambda \ge \kappa$. This
remark also applies to our following result that implies an
 analogue of theorem \ref{tm:limit} for $\mu$-extraresolvability
instead of $\mu$-resolvability.

\begin{theorem}\label{tm:meres1p}
Let ${\kappa}< {\lambda}=\cf({\lambda})\le (2^{\kappa})^+$  be  
infinite cardinals.
 Then there is a $\ccal(\kappa)$-space that is

\begin{enumerate}
\item
hereditarily ${\kappa}$-resolvable,
\item
 hereditarily ${\mu}$-\eres\ for all ${\mu}<{\lambda}$,
 \item
not ${\lambda}$-\eres.
\end{enumerate}
\end{theorem}

\begin{proof}
Similarly as in the proof of \ref{tm:limit}, let the sequence
$\vec {\lambda}=\<{\lambda}_{\zeta}:{\zeta}<{\lambda}\>$ be given
by $\lambda_\zeta = \omega_\zeta$
if ${\lambda}$ is a limit (hence inaccessible) cardinal, and  let
${\lambda}_{\zeta}={\rho}$ for each ${\zeta}<{\lambda}$ if
${\lambda}={\rho}^+$ is a successor.

Using fact \ref{tm:ld} again, we can find the following two types
of families of $2$-partitions of ${\kappa}$: $$\mathbb
B=\bigl\{\<B_{\xi}^0,B_{\xi}^1\>:{\xi}<2^{\kappa}\bigr\}$$ and
$$\dbb_{\zeta}=\bigl\{\<D_{{\zeta},\nu}^0,D_{{\zeta},\nu}^1\>:
{\nu}<{\lambda}_{\zeta}\bigr\}$$ for all ${\zeta}<{\lambda}$,
moreover a countable family $$\mathbb
G=\{\<G^i_n:i<{\kappa}\>:n<{\omega}\}$$ of ${\kappa}$-partitions
of ${\kappa}$
such that $\bbb\cup
\bigcup_{{\zeta}<{\lambda}}\dbb_{\zeta}\cup\mathbb G$ is
${\kappa}$-independent.

Now let $\dcal$ be the family of all sets of the form $\cap_{i<n}
E_i\cap \mathbb G[{\eta}]$ where $n < \omega$ and $E_i\in
\hint{\dbb_{{\zeta}_i}}$ with all the $\zeta_i$ distinct, moreover
${\eta}\in \fn({\omega},{\omega})$. It is easy to see that $\dcal$
is a family of ${\kappa}$-dense sets in $X_{\bbb}$, so we may
apply theorem \ref{tm:main} with $\bbb$ and $\dcal$ to get a
family of partitions $\cbb$ satisfying \ref{tm:main} (1) - (4). We
claim that $X_{\cbb}$ is as required.

Indeed, as we have already seen many times, the $\mathbb G[\eta]$
components of the elements of $\dcal$ can be used to show that
every $D\in\dcal$ is ${\kappa}$-resolvable. But then, as
$X_{\cbb}$ is both $\dcal$-forced and NODEC, every crowded
subspace of $X_{\cbb}$ is ${\kappa}$-resolvable by lemma 2.10,
hence (1) is proven.

To prove (2), we need the following statement.
\begin{claim}\label{cl:nwd4}
Assume that $\zeta < \lambda$ and $\{\nu,\nu'\} \in
[\lambda_\zeta]^2$. Then $$Y = D_{{\zeta},\nu}^0 \cap
D_{{\zeta},\nu'}^0\in \nwds(X_{\cbb}).$$
\end{claim}

\begin{proof}

Let  $D=\cap_{i<n} E_i\cap G$ be an arbitrary element of $\dcal$,
where $n\in {\omega}$, $\{{\zeta}_i:i<n\}\in \br {\lambda};n;$
with $\ E_i\in \hint{\dbb_{{\zeta}_i}}$ for all $i < n$, and $G =
\mathbb G[\eta]$ for some $\eta \in \fn(\omega,\omega)$. Our aim
is to check that $D\setm Y$ is dense, hence, by shrinking $D$ if
necessary, we may assume that ${\zeta}_0={\zeta}$ and
$E_0=D_{\zeta,\varphi}^0\setm \cup_{{\xi}\in\Psi}D_{\zeta,\xi}^0$.
Since ${\nu}\ne {\nu}'$ we can assume that ${\varphi}\ne {\nu}$.
Then
\begin{multline*}
D\setm Y\supset (\cap_{i<n} E_i\cap G)  \setminus
D_{{\zeta},\nu}^0 =\\=
(D_{{\zeta},\varphi}^0\setm
\cup_{{\xi}\in\Psi\cup\{{\nu}\}}D_{{\zeta},\xi}^0)
\cap\bigcap^{n-1}_{i=1}E_i\cap G\in\dcal.
\end{multline*}
Hence, $D\setm Y$ is indeed dense and so, by  lemma \ref{i:nwd},
$Y$ is nowhere dense in $X_{\cbb}$.
\end{proof}

Assume now that $D=\cap_{i<n} E_i\cap G$ is again an arbitrary
element of $\dcal$ with $\ E_i\in \hint{\dbb_{{\zeta}_i}}$ for all
$i < n$. By claim \ref{cl:nwd4}, for every $\zeta$ that is
distinct from all the $\zeta_i$ the collection $$\{D \cap
D_{{\zeta},\nu}^0 : \nu < \lambda_\zeta\}$$ consists of members of
$\dcal$ that have pairwise nowhere dense intersections, hence $D$
is ${\lambda}_{\zeta}$-\eres\ . Clearly, this implies that $D$
is ${\mu}$-\eres\ for all ${\mu}<{\lambda}$.  By lemma 2.10, since
$X_{\cbb}$ is $\dcal$-forced and NODEC it follows that $X_{\cbb}$
is hereditarily $\mu$-extraresolvable for all ${\mu}<{\lambda}$
and thus (2) has been established.

Finally, a standard $\Delta$-system and counting argument proves
that for each $\ecal\in \br\dcal;{\lambda};$ there is $\fcal\in
\br \ecal;{\lambda};$ such that $F\cap F'\in \dcal$ whenever
$\{F,F'\} \in \br \fcal;2;$. Hence, by lemma \ref{i:space}, the
space $X_{\cbb}$ is not ${\lambda}$-\eres , proving (3).
\end{proof}

Having seen these parallels between resolvability and
extraresolvability, it is interesting to note that we do not know
if the analogue of Bashkara Rao's ``compactness'' theorem  holds
for extraresolvability.

\begin{problem}
Assume that
 ${\lambda}$ is a singular cardinal with  $\cf({\lambda})={\omega}$ and
the space $X$ is
$\mu$-\eres\ for all $\mu<{\lambda}$. Is it true then that $X$ is
also ${\lambda}$-\eres\ ?
\end{problem}

Both theorems \ref{tm:meres1} and \ref{tm:meres1p} imply, in ZFC,
that for every infinite cardinal ${\kappa}$ there is a
$2^{\kappa}$-\eres\  $\ccal(\kappa)$-space.
However, theorem \ref{tm:meres2} below implies that this fails for
strong $2^\kappa$-extraresolvability. To prove \ref{tm:meres2},
however, we need some preparatory work.

\begin{definition}\label{df:hyres}
Let $\kappa$ be any cardinal. A topological space $X$ is called
{\em ${\kappa}$-\whyres} iff there is a $\kappa$-sequence $\langle
A_{\alpha}:{\alpha}<{\kappa}\rangle$ of pairwise disjoint elements
of $\br X;<{\kappa};$ such that $|A_{\alpha}|\le |{\alpha}|$ for
all $\alpha < \kappa$ and $\cup\{A_{\alpha}:{\alpha}\in I\}$ is
${\kappa}$-dense in $X$ whenever $I\in \br {\kappa};{\kappa};$.
If, in addition, $$\cup\{A_{\alpha}:{\alpha}\in K\} \in \ncal(X)$$
for each $K\in \br {\kappa};<{\kappa};$ then $X$ is called {\em
${\kappa}$-\hyres.} Finally, we say that $X$ is {\em \whyres\ 
(\hyres)} iff it is $\operatorname{\Delta}(X)$-\whyres\ 
($\operatorname{\Delta}(X)$-\hyres).
\end{definition}

We call a subfamily $\fcal$ of $\br {\kappa};{\kappa};$ {\em
boundedly} almost disjoint (BAD) if the intersection of any two
members of $\fcal$ is bounded in $\kappa$. Of course, if $\kappa$
is regular then any almost disjoint subfamily of $[\kappa]^\kappa$
is BAD. Moreover it is standard to show that for every infinite
$\kappa$ there is a BAD subfamily of $\br {\kappa};{\kappa};$ of
size ${\kappa}^+$. Thus from the above definitions we get the
following fact, explaining the term \hyres.

\begin{fact}
Any \hyres\ space $X$ is extraresolvable and if, in addition,
$\operatorname{\Delta}(X) = \nwd(X)$ then $X$ is \seres.
\end{fact}

\begin{definition}
Let $X$ be a topological space and ${\kappa}$ be an infinite
cardinal. A point $p\in X$ is said to be a {\em ${\kappa}$-\appr}
iff there is a one-to-one ${\kappa}$-sequence of points converging
to $p$ in $X$.
\end{definition}

\begin{lemma}\label{lm:appr}
If a topological space $X$ contains a  dense set $D$ of size $\le{\kappa}$
consisting of ${\kappa}$-\appr\ points then $X$ is ${\kappa}$-\whyres.
\end{lemma}

\begin{proof}\prlabel{lm:appr}
Enumerate $D$ as $\{d_{\zeta}:{\zeta}<{\kappa}\}$ with possible
repetitions. For each $d\in D$ fix a one-to-one sequence
$\{x_d({\xi}):{\xi}<{\kappa}\}\subs X\setm \{d\}$ converging to
$d$. By transfinite recursion on ${\alpha}<{\kappa}$ we may easily
construct a sequence $\<A_{\alpha}:{\alpha}<{\kappa}\>$ such that
\begin{enumerate}[(1)]
\item $A_{\alpha}\subs X\setm \bigcup\{A_{\delta}:{\delta}<{\alpha}\}$,
\item  $|A_{\alpha}|\le |{\alpha}|,$
\item  $A_{\alpha}\cap \{x_{d_{\zeta}}({\beta}):{\beta}\ge {\alpha}\}\ne \empt$
for all ${\zeta}\le{\alpha}$.
\end{enumerate}
It remains to show that $A_I=\bigcup\{A_{\alpha}:{\alpha}\in
I\}$ is $\kappa$-dense in $X$ whenever $I\in \br
{\kappa};{\kappa};$. So let $G\ne \emptyset$ be open and fix $d\in
D \cap G$. There is ${\zeta}<{\kappa}$ with $d_{{\zeta}}=d$. Then
for each ${\alpha}\in I\setm  {\zeta}$ we have $A_{\alpha}\cap
\{x_d({\beta}):{\beta}\ge {\alpha}\}\ne \empt$. But the sequence
$\{x_d({\xi}):{\xi}<{\kappa}\}$ is eventually in $G $ and the
$A_\alpha$'s are pairwise disjoint, consequently we have $|G \cap
 A_I| \ge \kappa.$
\end{proof}

The Cantor cube $\dis (2)^{2^{\kappa}}$ has a dense subset of
cardinality $\kappa$, moreover every point of $\dis
(2)^{2^{\kappa}}$ is a ${\kappa}$-limit point. Thus from lemma
\ref{lm:appr} we obtain the following fact.

\begin{fact}\label{tm:hyres}For each cardinal ${\kappa}$,
 the Cantor cube $\dis(2)^{2^{\kappa}}$ has a  ${\kappa}$-\whyres, dense subspace $X$
with $|X|=\operatorname{\Delta}(X)={\kappa}$.
\end{fact}

Using our main theorem \ref{tm:main} we can improve this as
follows.

\begin{theorem}\label{tm:d2hyres}
For each cardinal ${\kappa}$ there is a hyperresolvable  (and
hence strongly extraresolvable)  $\ccal(\kappa)$-space.
\end{theorem}

\begin{proof}\prlabel{tm:d2hyres}
By \ref{lm:embed} and  fact \ref{tm:hyres} we can find a
$\kappa$-independent family
$$\bbb = \{\langle B_{\xi}^0,B_{\xi}^1 \rangle : \xi <
2^{\kappa}\}$$ of $2$-partitions of ${\kappa}$ such that
$X_{\bbb}$ is ${\kappa}$-fragmented  by the sequence
$\acal=\<A_{\nu}:{\nu}<{\kappa}\>$.

As above, for any $I\subs {\kappa}$ we write
$A_I=\bigcup\{A_{\nu}:{\nu}\in I\}$. Then
\begin{displaymath}
  \dcal=\{A_I:I\in \br {\kappa};{\kappa};\}
\end{displaymath}
is a family of ${\kappa}$-dense sets in $X_{\bbb}$. So we can
apply theorem \ref{tm:main} to $\bbb$ and  $\dcal$ and get a
family $\cbb$ of $2$-partitions of $\kappa$ satisfying
\ref{tm:main} (1) - (4).

We claim that the sequence $\acal$ witnesses that $X_{\cbb}$ is
$(\operatorname{\Delta}(X_{\cbb}) =){\kappa}$-\hyres\  . Indeed,
every $A_I$ remains $\kappa$-dense in $X_{\cbb}$ for $I\in \br
{\kappa};{\kappa};$ because $A_I\in\dcal$. Moreover, if $J\in
\br{\kappa};<{\kappa};$ then for each $I\in
\br{\kappa};{\kappa};$, we have $A_I\setm A_J=A_{I\setm J}\in
\dcal$, consequently lemma \ref{i:nwd} may be applied to conclude
that
$A_J$ is nowhere dense in $X_{\ccal}$. 

\end{proof}

\begin{remark}
 The spaces obtained from  theorem
\ref{tm:d2hyres}  do not contain nontrivial convergent sequences
of any length because they are NODEC. This shows that the converse
of lemma \ref{lm:appr} fails.
\end{remark}

After this preparation we are now ready to formulate and prove
theorem \ref{tm:meres2}.

\begin{theorem}\label{tm:meres2}
Let ${\kappa} = cf(\kappa) < {\lambda}$ be two infinite cardinals.
Then the following three statements are equivalent:
\begin{enumerate}[(i)]
\item\label{lseres} There is a
\seresk{\lambda}\ but not $\lambda^+$-\eres\ 
$\ccal(\kappa)$-space.
\item\label{ttwo}  There is a \seresk{\lambda}\  space $X$ with $$|X|=\nwd(X)={\kappa}.$$

\item\label{ad} There is an almost disjoint family $\tcal\subs\br {\kappa};{\kappa};$
of size ${\lambda}$.
\end{enumerate}
\end{theorem}
\begin{proof}
\prtxtlabel{\ref{tm:meres2}} Clearly (i) implies (ii) implies
(iii). To prove that (iii) implies (i), we again use fact
\ref{tm:hyres} and observation \ref{lm:embed} to find an
independent, \separ\ family $\bbb = \{\langle B_{\xi}^0,B_{\xi}^1
\rangle : \xi < 2^{\kappa}\}$ of $2$-partitions of ${\kappa}$
such that $X_{\bbb}$ is ${\kappa}$-fragmented  by the sequence
$\acal=\<A_{\nu}:{\nu}<{\kappa}\>$.

Since  $A_I=\bigcup\{A_{\nu}:{\nu}\in I\}$ is $\kappa$-dense in
$X_{\cbb}$ for each $I \in [\kappa]^\kappa$, we may apply theorem
\ref{tm:main} to $\bbb$ and the family of $\kappa$-dense sets
$$\dcal=\{A_T : T\in \tcal\}$$
to get a family $\cbb$ of $2$-partitions of $\kappa$ that
satisfies  \ref{tm:main} (1) - (4).

Since $\kappa$ is regular, the family of dense sets $\{A_T:T\in
\tcal\} \subset [\kappa]^\kappa$ is almost disjoint because
$\tcal$ is. This, together with  $\nwd(X_{\cbb}) = \kappa$ clearly
implies that $X_{\cbb}$ is \seresk{\lambda}.

Moreover, as $|\dcal|={\lambda}$ and $\cel(X_{\cbb})={\omega}$,
lemma \ref{i:space} implies that $X_{\cbb}$ is not
${\lambda}^+$-\eres.

\end{proof}

It is known (see e. g. \cite{B}) that if one adds at least
${\omega}_3$ Cohen reals to a model of GCH then in the resulting
generic extension there is no almost disjoint subfamily of $ \br
\oo;\oo;$ of size $\omega_3$. Consequently, by theorem
\ref{tm:meres2}, in such a ZFC model, although $2^{\omega_1}$ is
as big as you wish, there is no \seresk{\omega_3}\ space $X$ with
$|X|=\nwd(X)=\oo$.

\end{document}